\documentclass[preprint,12pt]{imsart}

\RequirePackage[OT1]{fontenc}
\RequirePackage{amsthm,amsmath,natbib}
\RequirePackage[colorlinks,citecolor=blue,urlcolor=blue]{hyperref}

\usepackage{amssymb}
\usepackage{a4wide}
\usepackage[usenames,dvipsnames]{xcolor}
\usepackage{graphicx}
\usepackage{mathrsfs}
\usepackage{bm}
\usepackage{multicol}

\usepackage{mathrsfs}
\usepackage{verbatim}
\usepackage[colorlinks,citecolor=blue,urlcolor=blue]{hyperref}
\makeatletter
\newtheorem{lemma}{Lemma}[section]
\newtheorem{theorem}{Theorem}[section]

\newtheorem{remark}{Remark}

\numberwithin{equation}{section}
\allowdisplaybreaks       

\newcommand\gai[1]{{\color{black}#1}}

\newcommand\mgai[1]{{#1}}

\begin{document}

\def\di{D^{-1}}
\def\dij{D^{-1}_j}

\newcommand{\bqa}{\begin{eqnarray}}
  \newcommand{\eqa}{\end{eqnarray}}
\newcommand{\bqn}{\begin{equation}\begin{array}{lll}}
    \newcommand{\eqn}{\end{array}\end{equation}}
\newcommand{\um}{\underline{m}}
\newcommand{\be}{\begin{equation}}
  \newcommand{\ee}{\end{equation}}
\newcommand{\sjln}{{\sum\limits_{j=1}^{n}}}
\newcommand{\siln}{{\sum\limits_{i=1}^{n}}}
\newcommand{\bgma}{\bm{\gamma}}
\newcommand{\non}{\nonumber\\}
\newcommand{\bbA}{{\bf A}}
\newcommand{\bbC}{{\bf C}}
\newcommand{\bbD}{{\bf D}}
\newcommand{\rE}{{\rm E}}
\newcommand{\bbK}{{\bf K}}
\newcommand{\bbB}{{\bf B}}
\newcommand{\bbI}{{\bf I}}
\newcommand{\bbr}{{\bf r}}
\newcommand{\bbq}{{\bf q}}
\newcommand{\bbQ}{{\bf Q}}
\newcommand{\bbR}{{\bf R}}
\newcommand{\bbL}{{\bf L}}
\newcommand{\bbM}{{\bf M}}
\newcommand{\bbN}{{\bf N}}
\newcommand{\bbU}{{\bf U}}
\newcommand{\bbV}{{\bf V}}
\newcommand{\bgL}{\bm{\Lambda}}
\newcommand{\bbS}{{\bf S}}
\newcommand{\bbT}{{\bf T}}
\newcommand{\bbx}{{\bf x}}
\newcommand{\bbX}{{\bf X}}
\newcommand{\bby}{{\bf y}}
\newcommand{\rtr}{{\rm tr}}
\newcommand{\tr}{{\rm tr}}
\newcommand{\bA}{{\bf A}}
\newcommand{\bT}{{\bf T}}
\newcommand{\bbgS}{{\boldsymbol\bbT_{nr}}}
\newcommand\re[1]{{\color{red}#1}}
\newcommand{\cC}{{\mathcal{C}}}
\newcommand{\ep}{{\epsilon}}
\newcommand{\rCov}{{\rm Cov}}
\newcommand{\Cov}{{\rm Cov}}
\newcommand{\rVar}{{\rm Var}}
\newcommand{\darrow}{\downarrow}
\newcommand{\iparrow}{\stackrel{i.p.}{\rightarrow}}
\newcommand{\umn}{\underline{m}_n}
\newcommand{\diag}{{\rm diag}}
\newcommand{\rank}{{\rm rank}}
\newcommand{\rP}{{\rm P}}
\newcommand{\gs}{\Sigma}
\newcommand{\bbSigma}{{\boldsymbol{\Sigma}}}
\newcommand{\hSigma}{\hat{\bm{\bbT_{nr}}}}
\newcommand{\bSigma}{\bm{\bbT_{nr}}}
\renewcommand\hat[1]{\widehat{#1}}
\renewcommand\tilde[1]{\widetilde{#1}}

\newcommand\E{\mathbb{E}}
\newcommand\N{\mathcal{N}}
\newcommand\lb{\left(}
\newcommand\rb{\right)}
\newcommand\z{\mathbf{z}}
\newcommand\y{\mathbf{y}}
\newcommand\norm[1]{\left\lVert#1\right\rVert}
\renewcommand{\d}[1]{\ensuremath{\operatorname{d}\!{#1}}}
\newcommand{\rvv}[1]{{#1}}
\newcommand{\rdd}{{}} 
\newcommand{\rddd}{{}}    

\newcommand{\bgai}[1]{\color{red}#1}

\begin{frontmatter}
	\title{Joint CLT for eigenvalue  statistics from several dependent 
		large dimensional sample covariance matrices with application}
	\runtitle{Joint CLT for linear spectral statistics}
	
	\begin{aug}
      \author{\fnms{Weiming} \snm{Li}\thanksref{a,e1}\ead[label=e1,mark]{li.weiming@shufe.edu.cn} }
      \author{\fnms{Zeng} \snm{Li}\thanksref{c,e2}\ead[label=e2,mark]{zxl278@psu.edu}}
      \and
      \author{\fnms{Jianfeng} \snm{Yao}\thanksref{b,e3}\ead[label=e3,mark]{jeffyao@hku.hk}}

      \runauthor{Weiming Li, Zeng Li  and Jianfeng Yao}

      \affiliation{Shanghai University of Finance and Economics,
        Pennsylvania State University and The University of Hong Kong}

      \begin{multicols}{3}
        \address[a]{School of Statistics and Management\\
          Shanghai University of Finance and Economics\\
          \printead{e1}}

        \address[c]{Department of Statistics\\ Pennsylvania State University\\
          \printead{e2}}

        \address[b]{Department of Statistics and Actuarial Science\\
          The University of Hong Kong\\
          \printead{e3}}
      \end{multicols}      

	\end{aug}
	
	\begin{abstract}
      Let $\bbX_n=(x_{ij})$ be a $k \times n$ data matrix with complex-valued,
      independent and  standardized entries 
      satisfying a Lindeberg-type moment condition.  We consider
      simultaneously $R$ sample covariance matrices
      $\bbB_{nr}=\frac1n \bbQ_r \bbX_n \bbX_n^*\bbQ_r^\top,~1\le r\le
      R$, where the $\bbQ_{r}$'s  are \rdd{nonrandom} real
      matrices with common dimensions $p\times k~(k\geq p)$.
      Assuming that both the dimension $p$ and the sample size $n$ grow to
      infinity, the
      limiting distributions of the eigenvalues of the matrices  $\{\bbB_{nr}\}$
      are identified, and 
      as 
      the main result of the paper, we establish  a  joint central limit theorem  for linear
      spectral statistics of the $R$ matrices $\{\bbB_{nr}\}$.
      Next, this new CLT is applied  to the problem of testing a high dimensional white noise
      in time series modelling.
      \rddd{
        In experiments the derived test has a controlled size and is
        significantly faster than the classical permutation test,
        though it does have lower power.}
      This
      application highlights the necessity of such joint CLT 
      in  the presence of several dependent sample covariance
      matrices. In contrast,  
      all the existing \rdd{works} on CLT for linear spectral statistics
      of large sample covariance \rdd{matrices} deal with a single sample covariance matrix ($R=1$). 	 
    \end{abstract}
	
	\begin{keyword}
		\kwd{Large  sample covariance matrices}
		\kwd{central limit theorem}
		\kwd{linear spectral statistics}
		\kwd{white noise test}
		\kwd{high-dimensional times series}
	\end{keyword}
	
\end{frontmatter}

{\sl MSC 2010 Mathematics Subject Classifications}: Primary 62H10;  secondary 60B12, 60B20.


\section{Introduction}
\label{sec:intro}

Modern information technology tremendously accelerates computing speed and greatly enlarges the amount of data storage, which enables us to collect, store and analyze data of large dimensions. Classical limit theorems in multivariate analysis, which normally assume fixed dimensions, become no longer applicable for dealing with high dimensional problems. Random matrix theory investigates the spectral properties of random matrices when their dimensions tend to infinity and hence provides a powerful framework for solving high dimensional problems. This theory has made systematic corrections to many classical multivariate statistical procedures in the past decades, see the monographs of \cite{BS10},  \cite{Yao15} and the review papers \rdd{\cite{John07} and \cite{PA14}}. It has found \rdd{diverse} applications in various research areas, including signal processing, network security, image processing, \rdd{statistical genetics} and other financial econometrics problems.

The sample covariance matrix \gai{is of central importance} in
multivariate analysis. Many \gai{fundamental} statistics in
multivariate analysis can be written as functionals of eigenvalues of
a sample covariance matrix $\bbS_n$ such as  linear spectral statistics
\rdd{(LSSs)} of the form \rdd{$f(\lambda_1)+\cdots+f(\lambda_p)$}
where the $\lambda_j$'s are eigenvalues of $\bbS_n$ and
\rdd{$f(\cdot)$} is a smooth function. The wide range of creditable
applications in high dimensional statistics triggered an uptick in the
demand for CLTs of such LSSs. Actually one of the most widely used
results in this area is \citet{BS04}, which \rdd{considers} a sample covariance matrix
\gai{of the form} 
$\bbB_n=\frac 1n
\mathbf{T}^{1/2}\bbX_n\bbX_n^*\mathbf{T}^{1/2}$, \rdd{where $\bbX_n=(x_{ij})$ is $p\times n$ matrix consisting of i.i.d.\ complex standardized entries and $\mathbf{T}$ is a $p\times p$ nonnegative Hermitian matrix.}
A  CLT for \rdd{LSSs} of $\bbB_n$ is established 
under the so-called {\em Mar$\check{c}$enko-Pastur regime}, \rdd{i.e. $n, p\rightarrow \infty,~ p/n\rightarrow c>0$.} Further refinement and extensions can be found in \rdd{\cite{ZBY2015}, \cite{CP15}, \cite{ZBY-cltF}, and \cite{ZBYZ16}.} Among them, \cite{ZBY2015} studied the unbiased sample covariance matrix when the population mean vector is unknown. \citet{CP15} looked into the ultra-high dimensional case when the dimension $p$ is much larger than the sample size $n$, that is \rdd{$p/n\rightarrow \infty$ as $n\rightarrow \infty$.
\citet{ZBY-cltF} derived the CLT for \rdd{LSSs} of large dimensional general Fisher matrices.
\citet{ZBYZ16} attempted to \rdd{relax} the fourth order moment condition in \citet{BS04} and incorporated it into the limiting parameters. }

However, this rich literature all deals with a {\em single} sample
covariance matrix $\bbB_n$. 
This paper, on the contrary, aims at the joint limiting behaviour of functionals of several groups of eigenvalues coming from different yet closely related sample covariance matrices. Specifically, we consider data samples   $\{\bby_{jr}\}_{1\leq j\leq n, 1\le  r\le R}$ \mgai{of the form}  $\bby_{jr}=\bbQ_r\bbx_j$ where
\begin{enumerate}
	\item[(M1)] \rdd{$\{\bbx_j,~1\le j\le n\}$} is a
	sequence of \rdd{$k$-dimensional independent and complex-valued random
	vectors} with independent standardized components \rdd{$(x_{ij})$},
	i.e. $\rE x_{ij}=0$ and $\rE|x_{ij}|^2=1$, and the dimension $k\ge
	p$;
	\item[(M2)]  $\{\bbQ_r, 1\leq r\leq R\}$ are $R$ \rdd{nonrandom} real
      matrices with common dimensions $p\times k$. The $R$ population covariance matrices $\{\bbT_{nr}=\bbQ_r\bbQ_r^\top, r=1,\ldots R\}$ are assumed product-commutative.
\end{enumerate}
We thus consider $R$ sample covariance matrices given by
\begin{equation}
\label{eq:Bn}
\bbB_{nr}=\rdd{\frac{1}{n}}\sum_{j=1}^{n}\bby_{jr} \bby_{jr}^*=\rdd{\frac{1}{n}}\bbQ_r \bbX_n \bbX_n^*\bbQ_r^\top~,~1\leq r\leq R,
\end{equation}
where $\bbX_n=(\bbx_1,\ldots,\bbx_n)$ is of size  $k\times n $, \rdd{$*$ denotes the conjugate transpose of matrices or vectors, and $\top$ stands for the transpose of real ones}.   Let $\lb \lambda_{jr} \rb_{1\leq j\leq p}$ be
the eigenvalues of $\bbB_{nr}$ $(1\leq r\leq R)$,  and consider $L\times R$
real-valued functions $\lb f_{lr}\rb_{1\leq l\leq L,~1\leq
  r\leq R}$. This leads to the family of $L\times R$ \rdd{LSSs}
\[  \varphi_{lr}=
f_{lr}(\lambda_{1r})+\cdots+f_{lr}(\lambda_{pr}), ~1\leq l\leq
L,~1\leq r\leq R~.
\]  
This paper establishes a joint CLT for these $L\times R$ statistics
$\{\varphi_{lr}\}$ under appropriate conditions. 

The importance of such joint CLT for \rdd{\rdd{LSSs} is} best explained and illustrated with the
following problem of testing a high dimensional white noise. 
Indeed, our motivation for \rdd{the joint CLT originates} from this
application to high-dimensional time series analysis.
Testing for white noise is a classical yet important problem in statistics, especially for diagnostic checks in time series modelling. For high dimensional time series, \rdd{current literatures focus} on estimation and dimension-reduction aspects of the modelling, including high
dimensional VAR models and various factor models. Yet model
diagnostics \rdd{have} largely been untouched. Classical omnibus tests such
as the multivariate Hosking and Li-McLeod tests are no longer suitable
for handling high dimensional time series. They become extremely
conservative, losing size and power dramatically. In a very recent
work, \cite{Li16} looked into this high dimensional portmanteau test
problem and proposed several new test statistics based on
single-lagged and multi-lagged sample auto-covariance matrices. More
precisely, \rdd{let's consider a $p$-dimensional  time
series modelled as \rvv{a linear process}
\begin{equation}\label{tsm}
  \bbx_t=\sum_{l\geq0}\bbA_l{\bf z}_{t-l},
\end{equation}
where $\{{\bf z}_t\}$ is a sequence of independent \rdd{$p$-}dimensional
random vectors with independent components
$\z_t=(z_{it})$ satisfying $\rE z_{it}=0,~\rE|z_{it}|^2=1,~
\rE|z_{it}|^4<\infty$.
Hence $\{\bbx_t\}$ has $\rE\bbx_t = {\bf 0}$, and its
lag-$\tau$ auto-covariance matrix $\bbSigma_\tau= \Cov(\bbx_{t+\tau},\bbx_t)$ depends on $\tau$
only. In particular, $\bbSigma_0={\rm Var}(\bbx_t)$} denotes the population
covariance matrix of the series.
The goal is to test whether  $\bbx_t$ is a white noise, i.e.
\begin{equation}\label{eq:H0}
  H_0: ~~\Cov(\bbx_{t+\tau},\bbx_t)={\bf 0},~\tau=1,\cdots,q,
\end{equation}
where $q\geq 1$ is a prescribed constant integer. \rdd{Let $\bbx_1,\ldots,\bbx_n$ be a sample generated from the model \eqref{tsm}.  The {\em lag-$\tau$ sample
auto-covariance matrix} is}
\begin{equation}
  \label{eq:Sigma}
  \widehat{\bbSigma}_{\tau}=\frac
  {1}{n}\sum_{t=1}^{n} \bbx_t\bbx_{t-\tau}^*,
\end{equation}
\rdd{ where $\bbx_t=\bbx_{n+t}$ when $t\leq 0$.} \cite{Li16} proposed a test statistic based on $\widehat{\bbSigma}_\tau$. For any given constant integer $1\leq \tau\leq q$, \rdd{their} test statistic $\widetilde{L}_\tau$ was designed to test the specific lag-$\tau$ autocorrelation of the sequence, i.e.
\[
  \widetilde{L}_\tau=\sum_{j=1}^p
  \lambda_{j,\tau}^2=\mathrm{Tr}(\widetilde{\bbM}_\tau^*\widetilde{\bbM}_\tau),
\]
where $\{\lambda_{j,\tau},~j=1,\cdots,p\}$ are the eigenvalues of
\[\widetilde{\bbM}_{\tau}= \frac12
  \left(\widehat{\bbSigma}_{\tau} +\widehat{\bbSigma}_{\tau}^* \right)
  =\frac
  {1}{2n}\sum_{t=1}^{n}\lb\bbx_t\bbx_{t-\tau}^*+\bbx_{t-\tau}\bbx_t^*\rb,\]
  which is the {\em symmetrized lag-$\tau$ sample
  auto-covariance matrix}. 

Notice that in matrix form \rdd{$\widetilde{\bbM}_{\tau}=\frac {1}{2n}\bbX_n(\bbD_\tau+\bbD_\tau^\top)\bbX_n^*$,} where

\rddd{
  \[  \bbD_\tau =
    \begin{pmatrix}
      \mathbf{0}  & \bbI_{n-\tau}\\
      \bbI_{\tau} &  \mathbf{0} 
    \end{pmatrix}
  \]
where $\bbI_m$ denotes the $m$th order unit matrix.
}
They have proved that, under the null hypothesis, 
 in the simplest setting when $\bbx_t=\z_t$, the limiting distribution of the test statistic $\widetilde{L}_\tau$ is
\rdd{  \begin{equation*}
   \frac{n}{p}\widetilde{L}_\tau - \frac{p}{2} \xrightarrow{d} \N\lb \frac{1}{2},~1+\frac{3c(\nu_4-1)}{2}\rb.
  \end{equation*}}
  Here, $p,n\rightarrow \infty$ and $p/n\rightarrow c>0$ and $\nu_4=\rE|z_{it}|^4$. The null hypothesis should be rejected for large values of $\widetilde{L}_\tau$. Simulation results  also show that this test statistic is consistently more powerful than the Hosking and Li-McLeod tests even when the latter two have been size adjusted.

It can be seen that the test statistic $\widetilde{L}_\tau$ is \rdd{an} LSS of $\widetilde{\bbM}_\tau$, which can be studied with the CLT in \citet{BS04}. Indeed, the non-null eigenvalues of the sample covariance matrix $\bbS_n=\frac1n \bbT_p^{1/2}\bbX_n\bbX_n^*\bbT_p^{1/2}$ considered there are the same as the matrix $\underline{\bbS}_n=\frac1n\bbX_n\bbT_p\bbX_n^*$ which resembles to the matrix $\widetilde{\bbM}_\tau$. However, the test statistic $\widetilde{L}_\tau$ can only detect serial dependence in a single lag each time. In order to capture a multi-lag dependence structure, a naturally more effective way would be accumulating the lags and consider the statistic
\begin{equation}\label{eq:lq}
 \mathcal{L}_q=\sum_{\tau=1}^q\widetilde{L}_\tau=\sum_{\tau=1}^q\mathrm{Tr}(\widetilde{\bbM}_{\tau}\widetilde{\bbM}_{\tau}^*).
 \end{equation}
\gai{Note that the CLT in \cite{BS04} (or in its recent extensions)
  can only be used to study the correlations between different \rdd{LSSs} of
  a {\em given} $\widetilde{\bbM}_\tau$, while to derive the null
  distribution of $\mathcal{L}_q$, we need to study the correlations
  between \rdd{LSSs} of {\em several} covariance matrices
  $\widetilde{\bbM}_\tau, ~1\leq \tau\leq q$. Consequently, we need to
  resort to the joint CLT studied in this paper to characterize the
  correlations among $\{\widetilde{L}_\tau,~1\leq\tau\leq q\}$.}  
It is worth noticing that  \citet{Li16} proposed another multi-lagged
test statistic $U_q$ by stacking a number of consecutive observation
vectors. It will be shown in this paper that this test statistic
$U_q$ is essentially much less powerful 
than $\mathcal{L}_q$ considered here due to the data loss caused by
observation stacking. This superiority of $\mathcal{L}_q$ over $U_q$ 
demonstrates the necessity and significance of studying a joint CLT
for \rdd{LSSs} of several dependent sample covariance matrices as proposed in
this paper.

The rest of the paper is organized as follows. The main results of the joint CLT of \rdd{LSSs} of different sample covariance matrices are presented in Section \ref{sec:CLT}. \rdd{The} application on high dimensional white noise test is provided in Section \ref{sec:AppWN} to demonstrate the utility of this joint CLT. \rdd{Numerical} studies have also lent full support to the theoretical derivations. Technical lemmas and proofs are left to Section \ref{sec:proofs}.
\rvv{Finally, \texttt{\small Matlab} codes for reproducing simulations
  in the paper are available at: ~ \texttt{\small http://web.hku.hk/\~{}jeffyao/papersInfo.html}}.

\section{Joint CLT for linear spectral statistics of $\{\bbB_{nr}\}_{1\leq r\leq R}$}
\label{sec:CLT}

\subsection{Preliminary knowledge on LSDs of $\{\bbB_{nr}\}_{1\leq r\leq R}$}
\label{sec:LSD}

Recall that the {\em empirical spectral distribution} (ESD)
of a $p\times p$ square  matrix $\bbA$ is the probability
measure $F^\bbA=p^{-1} \sum_{i=1}^p\delta_{\lambda_i}$,
where the $\lambda_i$'s  are eigenvalues of $\bbA$ and $\delta_a$ denotes the
Dirac mass at point $a$. For any probability measure $F$ \mgai{on the real line},
its Stieltjes transform is defined by
\[
m(z)=\int\frac{1}{x-z} dF(x), \quad  z\in
\mathbb{C}^+,
\]
where $\mathbb{C}^+$  denotes the upper complex plane.

The assumptions needed for the existence of \rdd{{\em limiting spectral distributions} (LSDs)} of $\{\bbB_{nr}\}_{1\leq r\leq R}$ are as follows. The setup as well as the following \rdd{Lemma} \ref{thm1} are established in \rdd{\cite{ZBYZ16}}.
\begin{description}
\item[Assumption (a)] \rdd{Both dimensions $p$ and $n$ tend to infinity such  that $c_n=p/n\to c>0$ as $n\to\infty$.}
\item[Assumption (b)]
  Samples are $\{\bby_{jr}=\bbQ_r\bbx_j, j=1,\ldots,n, r=1,\ldots,R\}$,  where $\bbQ_r$ is $p\times k$, \rdd{$\bbx_j=(x_{1j},\ldots,x_{kj})^\top$ is $k\times 1$},
  and the dimension $k$ \rdd{($k\geq p$)} is arbitrary. Moreover,
  $\{x_{ij}, i=1,\ldots,k, j=1,\ldots,n\}$ is a $k\times n$ array of
  independent random variables, not necessarily
  identically distributed, with common moments
  $$
  \rdd{\rE} x_{ij}=0, ~~ \rdd{\rE}|x_{ij}^2|=1, 
  $$
  and satisfying the following Lindeberg-type condition:
  for each $\eta>0$,
  $$
  \frac1{pn\eta^2}\sum_{i=1}^k\sum_{j=1}^n \sum_{r=1}^R\|\bbq_{ir}\|^2\rE|x_{ij}^2|I\Big(|x_{ij}|>\eta \sqrt{n}/\|\bbq_{ir}\|\Big)\to0,
  $$
  where $\|\bbq_{ir}\|$ is the Euclidean norm of the $i$-th column vector
  $\bbq_{ir}$ of $\bbQ_r$.
\item[Assumption (c)] The ESD  $H_{nr}$ of the
  population covariance matrix $\bbT_{nr}=\bbQ_r\bbQ_r^\top$ converges weakly
  to a probability distribution $H_r$, $r=1,\ldots,R$.
 Also the sequence of the spectral norm of $(\bbT_{nr})$ is  bounded in $n$ and $r$.

\end{description}

\begin{lemma}\label{thm1}\rdd{[Theorem 2.1 of \cite{ZBYZ16}]}
  Under Assumptions \rdd{(a)-(c)}, almost surely, the  ESD
  $F_{nr}$ of $\bbB_{nr}$ weakly converges  to a \rdd{nonrandom}
  LSD  $F^{c,H_r}$. Moreover, \rdd{the Stieltjes transform $m_r(z)$ of $F^{c,H_r}$ is} the unique
  solution to the following Mar\v{c}enko-Pastur equation
  \begin{equation}
    \label{eq:MP}
    m_r(z)= \int\frac{1}{t[1-c-czm_r(z)]-z}dH_r(t)~,
  \end{equation}
  on the set
  $\{ m_r(z) \in\mathbb{C}:  -(1 - c)/z + cm_r(z)\in\mathbb{C}^+\}$.
\end{lemma}

Define the {\em companion LSD} of $\bbB_{nr}$ as
$$
\underline{F}^{c,H_r}=(1-c)\delta_0+cF^{c,H_r}.
$$
It is readily checked that
$ \underline{F}^{c,H_r}$ is the LSD of the {\em companion sample
  covariance matrix} $\underline{\bbB}_{nr}=n^{-1}\bbX_n^{*}\bbQ_r^\top\bbQ_r\bbX_n$ (which is $n\times n$), and its Stieltjes transform
$\underline{m}_r(z)=-(1-c)/z+cm_r(z)$ satisfies the so-called Silverstein equation
\begin{equation}
  z=-\frac1{\um_r(z)} +c\int\frac{t}{1+t\um_r(z)}dH_r(t).
  \label{eq:Silv}
\end{equation}

\subsection{Main Results}

Let \rdd{$\bbA$ and $\bbB$} be two real symmetric $p\times p$ matrices satisfying \rdd{$\bbA\bbB=\bbB\bbA$}. The two matrices can then be diagonalized simultaneously. 
We define the joint spectral distribution of \rdd{$(\bbA, \bbB)$} as the two-dimensional spectral distribution of the complex matrix \rdd{$\bbA+{\bf i}\bbB$}, i.e.,
$$
G(x,y)=\frac{1}{p}\#\{i\leq p, \Re(s_i)\leq x,\ \Im(s_i)\leq y\},
$$
where $(s_i)$ are the $p$ eigenvalues of \rdd{$\bbA+{\bf i}\bbB$} and $\# E$ denotes \rdd{the} cardinality of a set $E$.

Recall the random vector of $L\times R$ \rdd{LSSs} of $\bbB_{nr}$'s 
\begin{equation}\label{lss}
\left(\int f_{\ell r}(x)dF_{nr}(x) \right)_{1\leq \ell\leq L, 1\leq r\leq R},
\end{equation}
where $(F_{nr})$ are the corresponding empirical spectral distributions of $(\bbB_{nr})$ and $(f_{\ell r})$ are $L\times R$ measurable functions on the real line. Our aim in this section is to establish the joint distribution of \eqref{lss} under suitable conditions. The main results are presented as follows.

\begin{description}
\item[Assumption (d)]
  The variables $\{x_{ij}, i=1,\ldots,k, j=1,\ldots,n\}$ are independent, with common moments
  $$
  E x_{ij}=0, ~~ E|x_{ij}^2|=1, ~~
  \beta_x=E|x_{ij}^4|-|Ex_{ij}^2|^2-2, ~~\mbox{and}~~\alpha_x=|Ex_{ij}^2|^2,
  $$
  and satisfying the following Lindeberg-type condition: for each $\eta>0$
  {\be
    \frac1{pn\eta^6}\sum_{i=1}^k\sum_{j=1}^n\sum_{r=1}^R\|\bbq_{ir}\|^2\rE|x_{ij}^4|I\Big(|x_{ij}|>\eta\sqrt{n/\|\bbq_{ir}\|}\Big)\to0.
    \ee}

\item[Assumption (e)]
  {Either  $\beta_x=0$, or the mixing matrices $\{\bbQ_r\}$ are such that the matrices $\{\bbQ_r^\top\bbQ_r\}$ are diagonal  (therefore with arbitrary $\beta_x$).}

\item[Assumption (f)]  \rdd{The joint} spectral distribution $H_{nrs}$ of $\bbT_{nr}$ and $\bbT_{ns}$ converges weakly
  to a probability distribution $H_{rs}$, $1\leq r, s\leq R$.
\end{description}

\gai{The framework with Assumptions (d)-(e)-(f) is inspired by the one advocated in \citet{ZBYZ16}. However, \rdd{an} extension is necessary here since we are dealing with several random matrices simultaneously while only one matrix is considered in the reference.}

\begin{theorem} \label{thm2}
  Under Assumptions \rdd{(a)-(f)}, let  $f_{11},\ldots,f_{LR}$ be \rdd{$L\times R$} functions analytic on a complex domain containing
  \be
  [I_{(0<c<1)}(1-\sqrt{c})^2\liminf_n\lambda_{\min}^{\bbT},~~(1+\sqrt{c})^2\limsup_n\lambda_{\max}^{\bbT}]\label{tag1.4} \ee
  with $\bbT=\{\bbT_{nr}\}$,  and $\lambda_{\min}^{\bbT}$ and
  $\lambda_{\max}^{\bbT}$ \mgai{denoting the smallest and the
    largest  eigenvalue of all the matrices in $\bbT$, respectively}.
  Then, the random vector
  \begin{equation}\label{lss-1}
  p\left(\int f_{\ell r}(x)dF_{nr}(x)-\int f_{\ell r}(x)dF^{c_n, \rdd{H_{nr}}}(x)\right)_{1\leq \ell\leq L, 1\leq r\leq R}.
  \end{equation}
  converges to an \rdd{$(L\times R)$}-dimensional Gaussian random vector $(X_{f_{11}},\ldots,X_{f_{LR}})$. The mean function is
  \begin{eqnarray}
    \rE X_{f_{\ell r}}&=&-\frac{1}{2\pi {\bf i}}\oint_{{\cal C}_1}
    f_{\ell r}(z)g_1(z)\left[\frac{\alpha_x}{(1-g_2(z))(1-\alpha_xg_2(z))}+\frac{\beta_x}{1-g_2(z)}\right]dz,\nonumber
  \end{eqnarray}
 where
  \begin{eqnarray*}
g_1(z)=\int\frac{c\underline{m}_r^3(z)t^2}{(1+t\underline{m}_r(z))^{3}}dH_r(t)\quad\text{and}\quad g_2(z)=\int\frac{c\underline{m}_r^2(z)t^2}{(1+t\underline{m}_r(z))^{2}}dH_r(t).
  \end{eqnarray*}
 The covariance function is
  \begin{eqnarray}
    {\rm Cov}(X_{f_{\ell'r}}, X_{f_{\ell s}})
    &=&\frac{1}{4\pi^2}\oint\limits_{{\cal
        C}_1}\oint\limits_{{\cal
        C}_2} f_{\ell'}(z_1)f_{\ell}(z_2)\frac{\partial^2g(z_1,z_2)}{\partial z_1\partial z_2}dz_1dz_2,
  \end{eqnarray}
  where $g(z_1,z_2)=\log(1-a(z_1,z_2))+\log(1-\alpha_x\rdd{a(z_1,z_2)})-\beta_xa(z_1,z_2)$ with
  \begin{eqnarray*}
  a(z_1,z_2)=\iint\frac{c\um_r(z_1)\um_s(z_2)t_1t_2}{(1+t_1\um_r(z_1))(1+t_2\um_s(z_2))}dH_{rs}(t_1,t_2).
 \end{eqnarray*}
 The contours $\mathcal C_1$ and $\mathcal C_2$ are non-overlapping, closed, \rdd{positively} orientated in the complex plane, and enclosing both the supports of $F^{ c, H_r}$ and of $F^{c,H_s}$.
\end{theorem}
\rdd{
\begin{remark}\label{remark1}
The centralization term in \eqref{lss-1} is the expectation of $f$ with respect to the distribution $F^{c_n, H_{nr}}$. This distribution is a finite dimensional version of the LSD $F^{c, H_{r}}$, which is defined by \eqref{eq:MP} with the parameters $(c,H_r)$ replaced with $(c_n, H_{nr})$. The use of $F^{c_n, H_{nr}}$ instead of $F^{c, H_{r}}$ aims to eliminate the effect of the convergence rate of $(c_n, H_{nr})$ to $(c, H_r)$.   
\end{remark}}
\rdd{
As an illustrative example of Theorem \ref{thm2}, we consider a simplified case where only two sample covariance matrices are involved, i.e. $\bbX_n\bbX_n'/n$ and $\bbQ\bbX_n\bbX_n'\bbQ^\top/n$, where $\bbX_n$ is a $p\times n$ matrix of i.i.d.\ real standard Gaussian variables. The corresponding population covariance matrices are $\bbI_p$ and $\bbT_n:=\bbQ\bbQ^\top$, respectively. It's clear that the ESD and its limit of the identity matrix $\bbI_p$ are both the Dirac measure $\delta_1$. Those of $\bbT_n$ are general and denoted by $H_n$ and $H$, respectively. Moreover, the joint spectral distribution function $H_{n12}(t_1,t_2)$ of $\bbI_p$ and $\bbT_n$ is equal to $H_n(t_2)$ for $t_1=1$ and zero otherwise. Denote the ESDs of the two sample covariance matrices by $F_{n1}$ and $F_{n2}$, respectively, and let 
$$
G_{n1}(x)=F_{n1}(x)-F^{c_n, \delta_{1}}(x)\quad\text{and}\quad G_{n2}(x)=F_{n2}(x)-F^{c_n, H_{n2}}(x).
$$
Then for any analytic function $f$, we have
\begin{eqnarray}
p\left(\int f(x)dG_{n1}(x), \int f(x)dG_{n2}(x)\right)^\top
\xrightarrow{d}
\N\left(\left(
\begin{matrix}
	v_1\\
	v_2
\end{matrix}
\right),
\left(\begin{matrix}
	\psi_{11}&\psi_{12}\\
	\psi_{12}&\psi_{22}
\end{matrix}
\right)\right).\label{case}
\end{eqnarray}
The parameters $(v_1, v_2, \psi_{11}, \psi_{22})$ of the marginal distributions in \eqref{case} have been derived by many authors, see \cite{BS04} and \citet{ZBYZ16} for example. While the covariance parameter $\psi_{12}$ is new and, from Theorem \ref{thm2}, it can be formulated as
\begin{eqnarray*}
v_{12} & =& -\frac{c}{2\pi^{2}}\oint_{\mathcal{C}_{1}}\oint_{\mathcal{C}_{2}}\frac{f(z_1)f(z_2)(\um_1(z_1)+z_1\um_1'(z_1))(\um_2(z_2)+z_2\um_2'(z_2))}{\left[c-(1+z_1\um_1(z_1)(1+z_2\um_2(z_2)\right]^2}\mathrm{d}z_{1}\mathrm{d}z_{2},
\end{eqnarray*}
where $\um_1(z)$ and $\um_2(z)$ are the companion Stieltjes transforms of the LSDs $F^{c,\delta_1}(x)$ and $F^{c,H}(x)$, respectively, and $\um'(z)$ denotes the derivative of $\um(z)$ with respect to $z$. For the simplest function $f(z)=z$, one may figure out $v_{12}=2c\int tdH(t)$ by the residual theorem. 
}

\section{Application to high dimensional white noise test}\label{sec:AppWN}

As discussed in \rdd{the introduction}, a notable application of the joint CLT presented in this paper is \rdd{to the high dimensional}  white noise test. In particular, it is expected that testing power could be gained by 
accumulating information across different lags, that is, by using the test statistic 
$
\mathcal{L}_q=\sum_{\tau=1}^q\mathrm{Tr}(\widetilde{\bbM}_{\tau}\widetilde{\bbM}_{\tau}^*)$
defined in \eqref{eq:lq}. 

Define the scaled statistic
\begin{equation}
  \label{eq:phi-q}
  \phi_q  = \frac{n}{p}\mathcal{L}_q - \frac{qp}{2}.
\end{equation}
The null hypothesis will be rejected for large values of $\phi_q$. We consider  high-dimensional  situations where the dimension $p$
is large compared to the sample size $n$. 
By \rdd{applying} the CLT in Theorem \ref{thm2}, the asymptotic \rdd{null} distribution of $\phi_q$ is derived as follows.

\begin{theorem}\label{MainThmMult}
  Let $q\geq 1$ be a fixed integer, and assume that
  \begin{enumerate}
  \item
    $\{z_{it},~i=1,\cdots,p,~t=1,\cdots,n\}$ \rdd{is a set of i.i.d.\ real-valued variables satisfying
    $\rE z_{it}=0,~\rE z_{it}^2=1,~\rE z_{it}^4=\nu_4<\infty$;}
  \item 
    \rvv{Relaxed  Mar\v{c}enko-Pastur regime: both 
      the sample size $n$ and the dimension $p$ grow to infinity such
      that 
      \[     0  < \liminf_{n\to\infty} \frac{p}{n}    \le  \limsup_{n\to\infty} \frac{p}{n}<\infty.
      \]
    }
  \end{enumerate}
  \rvv{
  Then in the simplest setting \rdd{where} $\bbx_t=\z_t$, we have 
    \begin{equation}\label{Thm2DEV}
      s(c_n)^{-1/2}  \{ \phi_q  -\frac{q}{2} \}
      \xrightarrow{d} \N (0,1),
    \end{equation}
  where $s(u)= q+ u(\nu_4-1)(q^2+q/2).$
  }
\end{theorem}

The proof of this theorem is given in Section \ref{sec:proofs}.

\medskip
Let $Z_\alpha$ be the upper-$\alpha$ quantile of the standard normal
distribution at level $\alpha$.
Based on Theorem~\ref{MainThmMult}, we obtain a  procedure for  testing  the null hypothesis in
\eqref{eq:H0} as follows.
\begin{equation}
  \label{eq:cLq-test}
  \text{\em Multi-Lag-$q$ test:\quad  Reject $H_0$ ~~if ~~}
  \phi_q -\frac{q}{2} > Z_\alpha \rvv{\sqrt{s(c_n)}}.
\end{equation}

\subsection{Simulation Experiments}\label{ssec:simul}

Most of the experiments of this section are designed to compare our test procedure in \eqref{eq:cLq-test} and the procedure based on the test statistic \rdd{$U_q$ from \citet{Li16} using Simes' method \citep{Simes86}. In \cite{Li16},} several testing procedures are discussed and the \rdd{test $U_q$} performs quite satisfactorily in terms of both size and power across different scenarios.

More precisely, let $q\geq 1$ be a fixed integer, define \rdd{ $p(q+1)$-dimensional vectors}
$\y_j=\lb
 \begin{array}{c}
  \bbx_{j(q+1)-q}\\
  \vdots\\
  \bbx_{j(q+1)}
\end{array}\rb,$ \rdd{$j=1,\ldots,N$}, $N=\left[\frac{n}{q+1}\right]$. Since \rdd{$\rE\bbx_t=0$ and
$\bbSigma_\tau=\Cov(\bbx_{t+\tau},\bbx_t)$,} we have
$$\Cov(\y_j)=\lb
\begin{array}{cccc}
  \bbSigma_0 & \bbSigma_1 & \cdots & \bbSigma_q \\
  \bbSigma_1 & \bbSigma_0 & \ddots & \vdots \\
  \vdots & \ddots & \ddots & \bbSigma_1 \\
  \bbSigma_q & \cdots & \bbSigma_1 & \bbSigma_0
\end{array}
\rb_{(q+1)p\times(q+1)p}.$$
The null hypothesis \rdd{$H_0: ~\Cov(\bbx_{t+\tau},\bbx_t)={\bf 0},~\tau=1,\cdots,q$} becomes $H_0: \bbSigma_1=\cdots=\bbSigma_q={\bf 0}$, a test for a block diagonal covariance structure of the stacked sequence $\{\y_j\}$.

When $\bbSigma_0=\sigma^2I_p$, the white noise test of $\{\bbx_t\}$ reduces to a sphericity test of $\{\y_j\}$. The well known John's test statistic $U_q$ can be adopted for this purpose.
In our case, the corresponding John's test statistic $U_q$ is defined as
\[U_q=\cfrac{\frac{1}{p(q+1)}\sum_{i=1}^{p(q+1)}\lb l_{i,q}-\overline{l}_q\rb^2}{{\overline{l}_q}^2},\]
where $\{l_{i,q},~ i=1,\rdd{\ldots},p(q+1)\}$ are the eigenvalues of the sample covariance matrix $\bbS_q=\frac 1N\sum_{j=1}^N \y_j\y_j^*,$ and $\overline{l}_q$ is their average.

Notice however that the use of blocks above reduces the sample size $n$ to the number of blocks $N=\left[\frac{n}{q+1}\right]$. This may result in \rdd{a} certain loss of power for the test. To limit such loss of power, we adopt \rdd{Simes' method} for multiple hypothesis testing in \citet{Simes86}. To implement \rdd{Simes' method, we denote 
	$$
	\bby_j^{(0)}=\bby_j(\bbx_1,\ldots,\bbx_n)
	$$
as the previously defined stacked sample. Then we rotate the sample $(\bbx_j)$ and define a series of new stacked samples $\bby_j^{(k)}$ for $k=1,\ldots,q$, that is,
$$
\bby_j^{(k)}=\bby_j(\bbx_{k+1},\ldots,\bbx_n,\bbx_1,\ldots,\bbx_{k})
$$
Then John's test statistic $U_q$ can be calculated based on the $q+1$ samples, which results in $q+1$ different statistics $\{U_q^{(k)}\}$.}  Moreover, let $P_k, \; 0\leq k\leq q$\rdd{, denote} the (asymptotic)  P-value for the John's test with the $k$-th set of \rdd{$\y_j$'s, i.e.} $$P_{k}=1-\Phi\lb (NU_q^{(k)}-p(q+1)-\nu_4+2)/2\rb,$$ where $\Phi(\cdot)$ is the cumulative distribution function of the standard normal distribution.
  Let $P_{(1)}\leq \cdots\leq P_{(q+1)}$ be a permutation of $P_0,\ldots,P_q$. Then by the Simes method, we reject $H_0$ if $P_{(k)}\leq \frac{k}{q+1}\alpha$ at least for one $1\leq k\leq q+1$ for the nominal level $\alpha$.

To compare our test statistic $\phi_q$ with multi-lag-$q$ John's test statistic $U_q$, we set \rdd{the} significance level $\alpha=5\%$ and the critical regions of the two tests are
\begin{itemize}
  \item[(1)] Our test $\phi_q$: $\{\phi_q  > \frac{q}{2} + Z_{0.95}
    \sqrt{s(c_n)}  ~\}$;
      \item[(2)] Multi-lag-$q$ John's test $U_q$ (\rdd{using} Simes' method): $\{\mbox{at least for one } 1\leq k\leq q+1, P_{(k)}\leq \frac{k}{q+1}0.05\}$.
\end{itemize}
Data are generated following four different scenarios for comparison:

\begin{itemize}
  \item [(I)] Test size under Gaussian white noise:
  $\bbx_t=\z_t$, \rdd{$(\z_t)\stackrel{i.i.d.}{\sim} \N_p({\bf 0},{\bf I}_p)$;}
      \item[(II)]  Test size under Non-Gaussian white noise:
            $\bbx_t=\z_t-2$, \rdd{$(z_{it})\stackrel{i.i.d.}{\sim} \mbox{Gamma}(4,0.5)$, $\rE(z_{it})=2$, ${\rm Var}(z_{it})=1$,} $\nu_4(z_{it})=4.5$;
          \item[(III)] Test power under a Gaussian spherical AR(1) process:
                    \rdd{$\bbx_t=\bbA\bbx_{t-1}+\z_t$, $\bbA=a{\bf I}_p$, $a=0.1$, $(\z_t)\stackrel{i.i.d.}{\sim} \N_p({\bf 0},{\bf I}_p)$;}
          \item[(IV)] Test power under a Non-Gaussian spherical AR(1) process:
                    \rdd{$\bbx_t=\bbA\bbx_{t-1}+(\z_t-2)$, $\bbA=a{\bf I}_p$, $a=0.1$, $(z_{it})\stackrel{i.i.d.}{\sim} \mbox{Gamma}(4,0.5)$, $\rE(z_{it})=2$, Var$(z_{it})=1$, $\nu_4(z_{it})=4.5$.}
\end{itemize}

\begin{figure}[htbp]
  \includegraphics[width=1.0\textwidth,trim={0cm 7.8cm 0cm 7.8cm},clip]{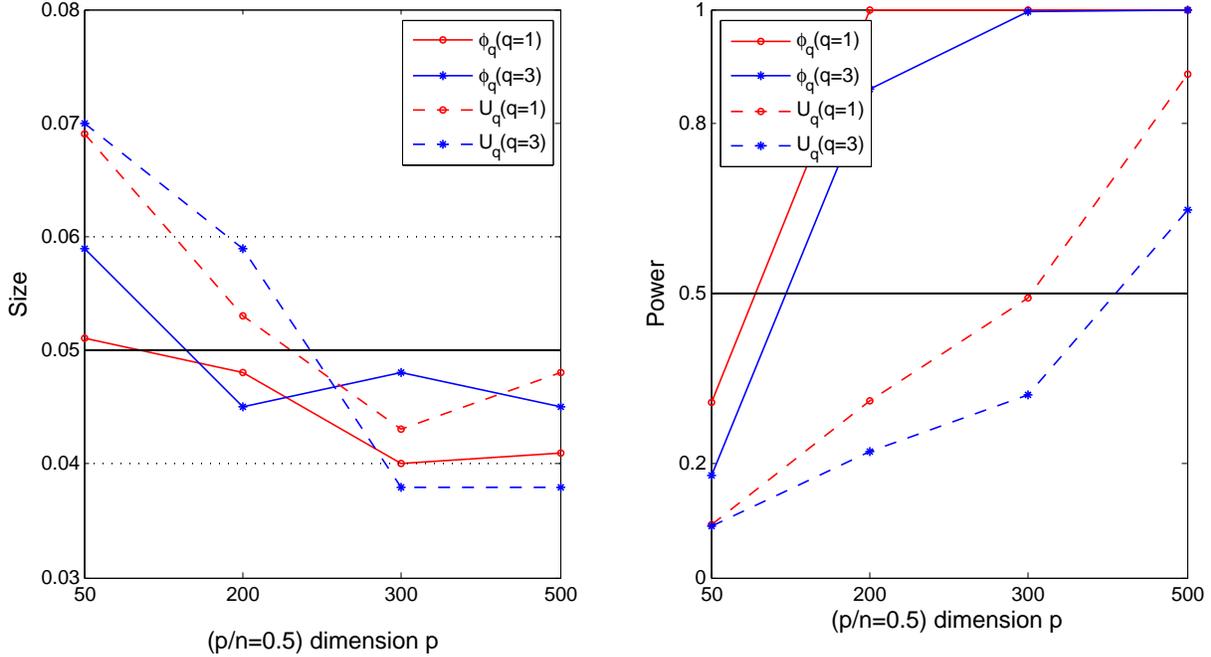}

  \caption{\small Sample plots of empirical sizes (left panel) from
    Table 1 (Scenario II with $c_n=0.5$) 
    and empirical powers (right panel) from Table 2 (Scenario III with
    $c_n=0.5$).}\label{fig:tableplot}
\end{figure}

\noindent
Various $(p,n)$-combinations are \rdd{tested to show the suitability
  of our test statistic for both low and high dimensional settings.}
Empirical statistics are obtained using 2000 independent
\rdd{replications}. Table \ref{Tab:LowSize} compares the empirical
sizes of \rdd{the two tests $\phi_q$ and $U_q$.}  It can be seen that
both \rdd{ of them} have reasonable sizes compared to the 5\% nominal
level across all the tested $(p,n)$-combinations. \rdd{Still, the} two
tests become slightly conservative under Non-Gaussian distributions in
Scenario (II) compared to \rdd{the} Gaussian case in Scenario (I). 
\rvv{A sample display of these sizes is given in Figure~\ref{fig:tableplot}
(left panel).}

In Table \ref{Tab:LowPower}, we compare the power of \rdd{the two tests.} Our \rdd{test} $\phi_q$ displays a generally much higher power than the multi-lag-$q$ John's test $U_q$, especially when \rdd{the dimensions $(p, n)$ become} larger. On the other hand, both tests have slightly lower power under the Non-Gaussian distribution than \rdd{under the} Gaussian distribution, which is consistent with the previous observation that \rdd{the} two tests become more conservative with Non-Gaussian populations.
\rvv{A sample display of these powers is given in Figure~\ref{fig:tableplot}
(right panel).}

To further explore \rdd{the powers} of the two \rdd{tests}, we varied the AR coefficient $a$ in Scenario (III) and (IV) from -0.1 to 0.1 ($a=0$ corresponds to \rdd{testing size}). Smaller values of the AR(1) coefficient $a$ are used here leading to a more difficult testing problem and a generally decreased power for both tests. Three dimensional settings are considered with $p/n\in\{ 0.1,0.5,1.5\}$ while the sample size is fixed as $n = 600$. \rdd{The number of independent replications is still 2000 in each case.} Results for Scenario (III) and (IV) are plotted in Figure 1. This Figure further consolidates that our \rdd{test $\phi_q$} dominates $U_q$ under all tested scenarios. A nonnegligible increase in the testing power of both test statistics as the \rdd{dimension $p$} becomes larger sheds more light on the blessings of high dimensionality. Still both tests are more conservative with Non-Gaussian population distribution than with Gaussian distribution.

\begin{figure}[htbp]
  \centering
  \includegraphics[width=1.0\textwidth,trim={0.5cm 3.5cm 0.5cm 3.5cm},clip]{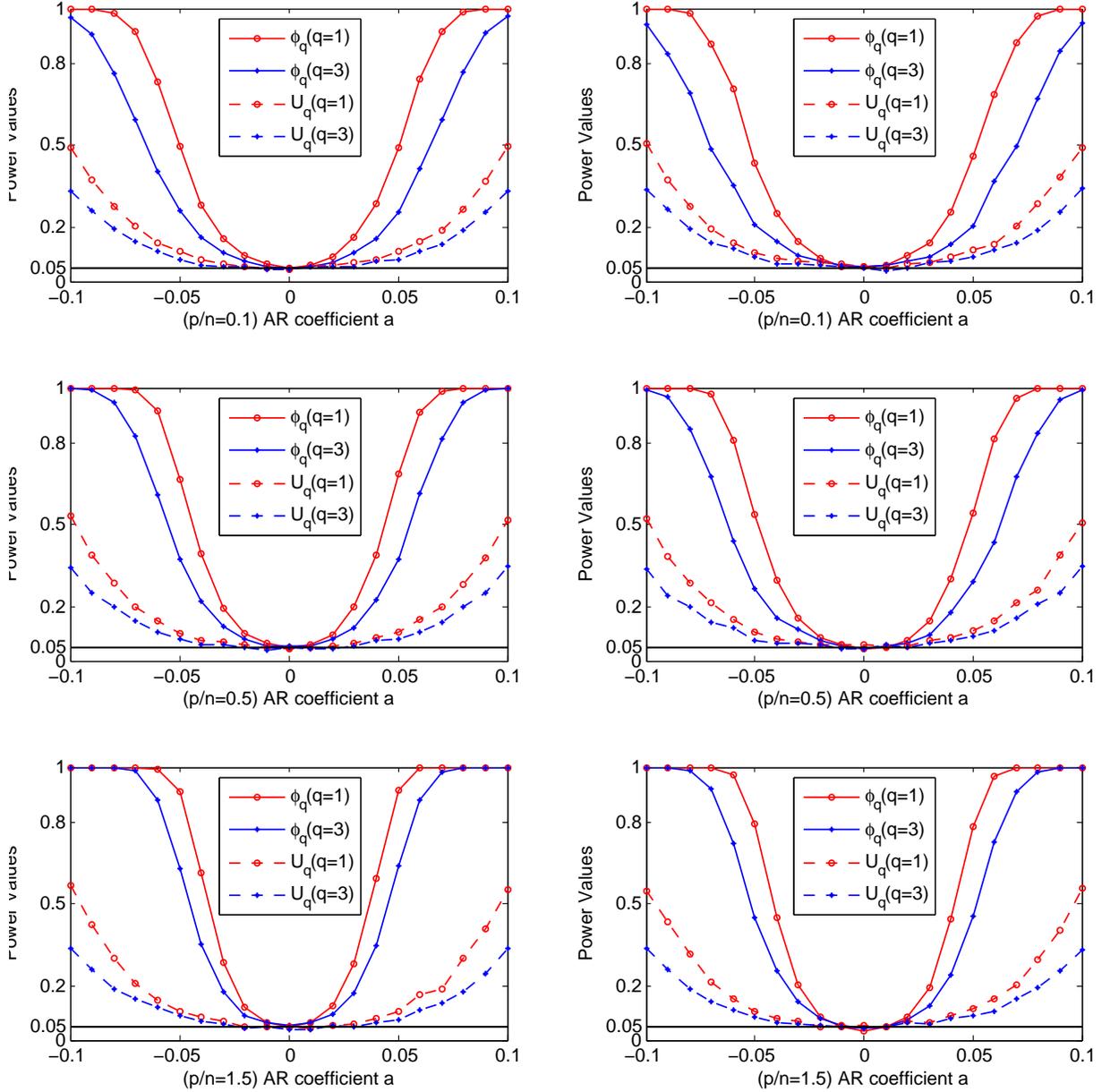}\\
  \caption{\small {Empirical \rdd{Powers for the two tests} with varying AR coefficient $a$ from -0.1 to 0.1. Left
      panel: Scenario(III) for Gaussian Distribution.  Right panel: Scenario(IV) for Gamma Distribution}\label{fig:powerline}}
\end{figure}

\subsection{Comparison to a permutation test}

As many complex analytic tools are employed to derive the asymptotic
null distributions of the test statistic $\phi_q$,  it is natural to
wonder about the performance of a ``simple-minded'' test procedure,
namely the permutation test. Under the null hypothesis of white noise,
since the sample vectors $\bbx_1,\ldots,\bbx_n$ have an  i.i.d. structure, one can permute these $n$
sample vectors say $B$ times to obtain an empirical upper 5\%
quantiles of the test statistic $\phi_q$. The null hypothesis will be
rejected if the  observed statistic $\phi_q$ from the original (non
permuted) sample vectors $\bbx_1,\ldots,\bbx_n$ is larger than this
empirical quantile.

Data are generated following the spherical AR(1) process in Scenario (III) and (IV) to compare this straightforward test with our test statistic $\phi_q$. In order to compare the power performance of two tests, the AR coefficient $a$ takes different values, $a=[0,0.05,0.09,0.1]$ ($a=0$ corresponds to testing size). The sample size is fixed as $n=300$ yet data dimension $p$ varies. As for the permutation test, the permutation times is set as $B=500$. The nominal level is $\alpha=5\%$. Testing size and power of two tests are shown in Tables \ref{tab:compareGauss} and \ref{tab:compareNonGauss} based on 500 replicates for all $(p,n)$ configurations.

\rddd{
  It can be seen that the  sizes of both tests are well controlled.
  As for their power, our test offers an acceptable
  level while   the permutation test consistently performs better
  in the tested cases.  
  However, the permutation test is extremely time consuming compared to our test.
  For instance,  to run one set of $(p,n)=(150,300)$ combination for
  500 replicates, it takes only {\bf 25 seconds} with our test, while
  almost {\bf 3 hours} for the permutation test with permutation times
  $B=500$. Particularly the computation time increases greatly when
  the sample size $n$ grows.
  Therefore,  our test statistic $\phi_q$  provides a very
  competitive choice for testing high dimensional white
  noise while the classical permutation test is simpler, more powerful
  though much slower.
}

\section{Proofs of the main theorems}\label{sec:proofs}

\subsection{Proof of Theorem \ref{thm2}}

\gai{The general strategy for our main Theorem \ref{thm2} follows the methods advocated in \citet{BS04}, with its most recent update in \citet{ZBYZ16}. However, as we are dealing with several random matrices simultaneously, all the technical steps for the implementation of this strategy have to be carefully rewritten. They are presented in this section.}

\subsubsection{\mgai{Sketch} of the proof of Theorem \ref{thm2}}

Let $v_0 > 0$ be arbitrary,  $x_r$ be any number greater than the right end point of interval \eqref{tag1.4}, and $x_l$ be any negative number if the left end point of \eqref{tag1.4} is zero, otherwise choose
$x_l\in(0,\liminf_{p\rightarrow\infty}\lambda_{\min}^{\bT}(1-\sqrt{c})^2)$. Define a contour $\mathcal C$ as
\begin{equation}\label{cont}
\mathcal C=\{x+iv: x\in \{x_r,x_l\}, v\in[-v_0,v_0]\}\cup \mathcal C_u \ \text{with}\ \mathcal C_u=\{x\pm iv_0: x\in[x_l,x_r]\},
\end{equation}
and let $\mathcal C_n=\mathcal C_u\cup \{x\pm iv: x\in\{x_l, x_r\}, v\in[n^{-1}\varepsilon_n, v_0]\}$ with
$\varepsilon_n\geq n^{-\alpha}$ for some $\alpha\in(0,1)$.
\rdd{By definition, the contour $\mathcal C$ encloses a rectangular region in the complex plane, which contains the union of the support sets of all the LSDs $F^{c, H_r}$, $1\leq r\leq R$. As a regularized version of $\mathcal C$, $\mathcal C_n$ excludes a small segment near the real line. }

\rdd{
	Let $m_{nr}(z)$, $\um_{nr}(z)$, $m_{nr}^0(z)$, $\um_{nr}^0(z)$ be the Stieltjes transforms of $F_{nr}$, $\underline{F}_{nr}$, $F^{c_n, H_{nr}}$, and $\underline{F}^{c_n, H_{nr}}$, respectively, where $F_{nr}$ is the ESD of $\bbB_{nr}$, $F^{c_n, H_{nr}}$ is the LSD defined in Remark \ref{remark1}, $F$ and $\underline{F}$ are linked by the equation $\underline{F}=(1-c_n)\delta_0+c_nF$.
    A major task of proving Theorem \ref{thm2} is to study the convergence of the empirical process
    $$M_{nr}(z):=p[m_{nr}(z)-m_{nr}^0(z)]=n[\um_{nr}(z)-\um_{nr}^0(z)].$$
    To this end, we need to truncate $M_{nr}(z)$ as
    \begin{eqnarray*}
    	\widehat M_{nr}(z)=
    	\begin{cases}
    		M_{nr}(z)&z\in \mathcal C_n,\\
    		M_{nr}(x+in^{-1}\varepsilon_n)& x\in \{x_l,x_r\}\quad\text{and}\quad v\in[0, n^{-1}\varepsilon_n],\\
    		M_{nr}(x-in^{-1}\varepsilon_n)& x\in \{x_l,x_r\}\quad\text{and}\quad v\in[-n^{-1}\varepsilon_n, 0],
    	\end{cases}
    \end{eqnarray*}
which agrees to $M_{nr}(z)$ on $\mathcal C_n$. This truncation is essential when proving the tightness $\widehat M_{nr}(z)$ on $\mathcal C$.
Write
$$
\widehat{\bbM}_{n}(z)=\left(\widehat{M}_{n1}(z),\ldots, \widehat{M}_{nR}(z)\right),
$$
we will establish its convergence as stated in the following lemma.}
\begin{lemma}\label{lem1t}
	Under Assumptions (a)-(f), 
	\rdd{$\widehat \bbM_n(\cdot)$} converges weakly to a Gaussian
	process $\rdd{\bbM\rdd{(\cdot)}}=(M_1,\ldots,M_R)(\cdot)$ on $\cC$. The mean function is
	\begin{eqnarray}
	\rE M_r(z)=\frac{\alpha_xg_1(z)}{(1-g_2(z))(1-\alpha_xg_2(z))}+\frac{\beta_xg_1(z)}{1-g_2(z)},\nonumber
	\end{eqnarray}
	where
	\begin{eqnarray*}
		g_1(z)=\int\frac{c\underline{m}_r^3(z)t^2}{(1+t\underline{m}_r(z))^{3}}dH_r(t)\quad\text{and}\quad g_2(z)=\int\frac{c\underline{m}_r^2(z)t^2}{(1+t\underline{m}_r(z))^{2}}dH_r(t).
	\end{eqnarray*}
	The covariance function is
	\begin{eqnarray*}
		{\rm Cov}(M_r(z_1), M_s(z_2))=
		-\frac{\partial^2}{\partial z_1\partial z_2}\left[\log(1-a(z_1,z_2))+\log(1-\alpha_x \rdd{a(z_1,z_2)})-\beta_x a(z_1,z_2)\right]
	\end{eqnarray*}
	where
\rdd{	\begin{eqnarray*}
		a(z_1,z_2)=c\iint\frac{t_1t_2\um_r(z_1)\um_s(z_2)}{(1+t_1\um_r(z_1))(1+t_2\um_s(z_2))}dH_{rs}(t_1,t_2).
	\end{eqnarray*}}
\end{lemma}
\rdd{From this lemma, Theorem \ref{thm2} follows by similar arguments \rdd{on} Pages 562 and 563 in \cite{BS04}.}

\subsubsection{Proof of Lemma \ref{lem1t}}
\label{sec:proofCLT}

\rdd{Following closely the steps of truncation, centralization and rescaling in Appendix B of \cite{ZBYZ16}, one may find that it is sufficient to prove this lemma under the assumption that 
\begin{equation}\label{xij}
|x_{ij}|<\frac{\eta_n\sqrt{n}}{\max_{1\leq r\leq R}\{\|\bbq_{ir}\|\}},
\end{equation}
where the constant $\eta_n\to 0$ as $n\to\infty$.}

Write for $r\in \{1,\ldots,R\}$ and 
$z\in\cC _n$, 
\rdd{\begin{eqnarray*}
M_{nr}(z)&=&p[m_{nr}(z)-\rE m_{nr}(z)]+p[\rE m_{nr}(z)-m_{nr}^0(z)]\\
&:=&M_{nr}^1(z)+M_{nr}^2(z).
\end{eqnarray*}}
The Lemma can be proved by verifying three conditions \citep{BS04}:
\begin{itemize}
	\item[] Condition 1: Finite dimensional convergence of $M_{n}^{1}(z)$ in distribution;
	\item[] Condition 2: Tightness of $M_{n}^{1}(z)$ on $\mathcal C_{n}$;
	\item[] Condition 3: Convergence of $M_{n}^{2}(z)$.  
\end{itemize}
Since the second and third conditions can be obtained directly from Lemma 5.1 in \cite{ZBYZ16}, we only consider the first one by showing that, for any \rdd{$W\times R$} complex numbers
$z_{11},\ldots, z_{WR}$, the random vector $(M_{nr}^{1}(z_{jr}))_{1\leq j\leq W, 1\leq r\leq R}$
converges to a Gaussian vector. Without loss of generality, we assume $\max\{\|\bbQ_r\|\}\leq1$. \rdd{We will also denote by $K$ any constants appearing in inequalities and $K$ may take on different values for different expressions.}

With the notation $\bbr_{jr}=(1/\sqrt{n})\bbQ_r\bbx_j$, we define some \rdd{quantities}:
\begin{align}
&\bbD_r(z)=\bbB_{nr}-z\bbI,\ \bbD_{jr}(z)=\bbD_r(z)-\bbr_{jr}\bbr_{jr}^*,\ \bbD_{ijr}(z)=\bbD_r(z)-\bbr_{ir}\bbr_{ir}^*-\bbr_{jr}\bbr_{jr}^*,\nonumber\\
&\ep_{jr}(z)=\bbr_{jr}^*\bbD^{-1}_{jr}(z)\bbr_{jr}-n^{-1}\rtr\bbT_{nr}\bbD^{-1}_{jr}(z),\nonumber\\
&\delta_{jr}(z)=\bbr_{jr}^*\bbD^{-2}_{jr}(z)\bbr_{jr}-n^{-1}\rtr\bbT_{nr}\bbD^{-2}_{jr}(z),\nonumber\\
&\beta_{jr}(z)=\frac{1}{1+\bbr_{jr}^*\bbD_{jr}^{-1}(z)\bbr_{jr}},\ \beta_{ijr}(z)=\frac{1}{1+\bbr_{ir}^*\bbD_{ijr}^{-1}(z)\bbr_{ir}},\ \bar{\beta}_{jr}(z)=\frac{1}{1+n^{-1}\rtr  \bbT_{nr}\bbD_{jr}^{-1}(z)},\nonumber\\
&b_{nr}(z)=\frac{1}{1+n^{-1}\rE\rtr  \bbT_{nr}\bbD_r^{-1}(z)},\ b_{12r}(z)=\frac{1}{1+n^{-1}\rE\tr \bbT_{nr} \bbD_{12r}^{-1}(z)},\nonumber
\end{align}
which will be frequently used in the sequel. \rdd{Note} that  quantities in the last two rows are all bounded in absolute value by $|z|/\Im(z)$.

By martingale difference decomposition, the process $M_{nr}^1(z)$ can be expressed as
\rdd{\begin{align*}
p[m_{nr}(z)-\rE m_{nr}(z)]
=&~\rtr [\bbD_r^{-1}(z)-\rE \bbD_r^{-1}(z)]\\
=&~\sum_{j=1}^n\rtr \rE_j[\bbD_r^{-1}(z)-\bbD^{-1}_{jr}(z)]-
\rtr \rE_{j-1}[\bbD_r^{-1}(z)-\bbD^{-1}_{jr}(z)]\\
=&~-\sum_{j=1}^n(\rE_j-\rE_{j-1})
\beta_{jr}(z)\bbr_{jr}^*\bbD^{-2}_{jr}(z)\bbr_{jr}\\
=&~-\frac{d}{dz}\sjln (\rE_j-\rE_{j-1})\log \beta_{jr}(z)\\
=&~\frac{d}{dz}\sjln (\rE_j-\rE_{j-1})\log \big(1+\ep_{jr}(z)\bar\beta_{jr}(z)\big),
\end{align*}
where the third equality is from the identity} $\bbD^{-1}_r(z)=\bbD_{jr}^{-1}(z)-\bbD_{jr}^{-1}(z)\bbr_{jr}\bbr_{jr}^*\bbD_{jr}^{-1}(z)\beta_{jr}(z)$ and \rdd{the last one is obtained using the identity $\beta_{jr}(z)=\bar{\beta}_{jr}(z)[1+\bar{\beta}_{jr}(z)\ep_{jr}(z)]^{-1}$.
We next show that 
\begin{equation*}
\frac{d}{dz}\sjln (\rE_j-\rE_{j-1})\log \big(1+\ep_{jr}(z)\bar\beta_{jr}(z)\big)-\frac{d}{dz}\sjln (\rE_j-\rE_{j-1})\ep_{jr}(z)\bar\beta_{jr}(z)=o_p(1).
\end{equation*}
Considering the second moment of the above difference, by the Cauchy integral formula, one may get
\begin{align}
~&~\rE\left|\frac{d}{dz}\sjln (\rE_j-\rE_{j-1})\big[\log \big(1+\ep_{jr}(z)\bar\beta_{jr}(z)\big)-\ep_{jr}(z)\bar\beta_{jr}(z)\big]\right|^2\non
=&~\rE\left|\sjln \frac1{2\pi {\bf i}}\oint_{|\zeta-z|=v/2}\frac{\big[\log \big(1+\ep_{jr}(\zeta)\bar\beta_{jr}(\zeta)\big)-\ep_{jr}(\zeta)\bar\beta_{jr}(\zeta)\big]}{(z-\zeta)^2}d\zeta\right|^2\non
\le&~\frac{K}{\pi^2 v^4}\sjln \oint_{|\zeta-z|=v/2}\rE|\ep_{jr}(\zeta)\bar\beta_{jr}(\zeta)|^4|d\zeta|.\label{eq:1}
\end{align}
From Lemma \ref{lem2} and the truncation in \eqref{xij}, we have
\begin{align*}
\rE|\ep_{jr}(\zeta)|^4\le &~\frac{K}{n^4}\bigg\{\rE\left[\rtr\bbT_{nr}\bbD_{jr}^{-1}(\zeta)\bbT_{nr}\bbD_{jr}^{-1}(\bar \zeta)\right]^2+\sum_{i=1}^k \rE|x_{ij}|^8 \rE|\bbq_{ir}^\top\bbD_{jr}^{-1}(\zeta)\bbq_{ir}|^4\bigg\}\non
\le&~ Kn^{-2}+K\eta_n^4n^{-1},
\end{align*}
by which the right hand side of \eqref{eq:1} tends to zero.
} 
Therefore, we need only to consider the limiting distribution of
\be
\frac{d}{dz}\sjln (\rE_j-\rE_{j-1})\ep_{jr}(z)\bar\beta_{jr}(z)=\frac{d}{dz}\sjln \rE_j\ep_{jr}(z)\bar\beta_{jr}(z)
\label{eqgclt}
\ee
in finite dimensional situations. 
\rdd{To verify the Lyapunov condition,  one can show that}
\begin{eqnarray}
&&\sjln\rE\left|(\rE_j-\rE_{j-1})\frac{d}{dz}\ep_{jr}(z)\bar\beta_{jr}(z)\right|^2 I\left(\left|(\rE_j-\rE_{j-1})\frac{d}{dz}\ep_{jr}(z)\bar\beta_{jr}(z)\right|\geq\epsilon\right)\non
&\leq&\frac{1}{\epsilon^2}\sjln \rE\Big|\rE_j\frac{d}{dz}\ep_{jr}(z)\bar\beta_{jr}(z)\Big|^4\to 0,\nonumber
\label{eqgclt2}
\end{eqnarray}
\rdd{where the convergence is again from Lemma \ref{lem2} and \eqref{xij}.
}
Hence, from the martingale \citep[Theorem 35.12]{Billingsley1995},
the random vector $(M_{nr}^1(z_{jr}))$ will tend to a Gaussian vector
$(M_r(z_{jr}))$ with covariance function
\begin{equation}
\rCov(M_r(z_1),M_s(z_2))=\lim_{n\to\infty}\sjln \rE_{j-1}\left(\rE_j\frac{\partial}{\partial z_1}\ep_{jr}(z_1)\bar\beta_{jr}(z_1)\cdot
\rE_j\frac{\partial}{\partial z_2}\ep_{js}(z_2)\bar\beta_{js}(z_2)\right).
\label{limcov}
\end{equation}
\rddd{We note that the referenced martingale CLT applies also to
  multidimensional martingale by considering arbitrary 
  linear combination of its components.}

Using the same approach of \cite{BS04} \rdd{on Page 571, one} may replace $\bar\beta_{jr}(z)$ by $b_{nr}(z)$. \rdd{Then}, by (1.15) of \cite{BS04}, we have
\begin{align}
\Gamma_{nrs}(z_1,z_2)\rdd{:=}&\sum_{j=1}^nb_{nr}(z_1)b_{ns}(z_2)\rE_{j-1}[\rE_j\epsilon_{jr}(z_1)
\rE_j\epsilon_{js}(z_2)]\non
=&~\frac1{n^2}\sum_{j=1}^nb_{nr}(z_1)b_{ns}(z_2)\Bigg[
\rtr\rE_{j}\bbQ_r^\top\bbD_{jr}^{-1}(z_1)\bbQ_r\rE_j\bbQ_s^\top\rdd{\bbD_{js}^{-1}}(z_2)\bbQ_s\non
&\quad\quad\quad\quad\quad\quad\quad\quad\quad
+\alpha_x\rtr\rE_{j}\bbQ_r^\top\bbD_{jr}^{-1}(z_1)\bbQ_r\rE_j\bbQ_s^\top\rdd{\bbD_{js}^{-1}}(z_2)\bbQ_s\non
&\quad\quad\quad\quad\quad\quad\quad\quad\quad
+\beta_x\sum_{i=1}^k
\bbq_{ir}^\top\rE_j\bbD_{jr}^{-1}(z_1)\bbq_{ir}\bbq_{is}^\top\rE_j\bbD_{js}^{-1}(z_2)\bbq_{is}\Bigg]\nonumber\\
&:=\Gamma_1+\alpha_x\Gamma_2+\beta_x\Gamma_3,
\label{bai2.6}
\end{align}
where $\alpha_x=|Ex_{11}^2|^2$ and $\beta_x=E|x_{11}^4|-\rdd{|Ex_{11}^2|^2}-2$.

Now we derive the limit of the first term in \eqref{bai2.6}. \rdd{The means is to replace $\bbD_{jr}^{-1}(z)$ (and similarly $\bbD_{js}^{-1}(z)$) by a proper nonrandom matrix.
For this, we introduce such a one}
\begin{eqnarray*}
	\bbL_{r}(z)=zI-\frac{n-1}{n}b_{12r}(z)\bbT_{nr},
\end{eqnarray*}
\rdd{whose inverse spectral norm is bounded, that is,}
\begin{equation}\label{z-b}
||\bbL_r(z)||^{-1}\leq \frac{|b_{12r}^{-1}(z)|}{\Im(zb_{12r}^{-1}(z))}\leq\frac{1+p/(nv)}{v}.
\end{equation}
\rdd{We will show that the major part of $\bbD_{jr}^{-1}(z)$ is just $-\bbL_{r}^{-1}(z)$.} From the identity $\bbr_{ir}^*\bbD_{jr}^{-1}(z)=\beta_{ijr}(z)\bbr_{ir}^*\bbD_{ijr}^{-1}(z)$, we get
\begin{align}
\bbD_{jr}^{-1}(z)+\bbL_r^{-1}(z)=&~\bbL_r^{-1}(z)\left(\bbD_{jr}(z)+\bbL_r(z)\right)\bbD_{jr}^{-1}(z)\nonumber\\
=&~\bbL_r^{-1}(z)\left(\sum_{i\neq j}\bbr_{ir}\bbr_{ir}^*-\frac{n-1}{n}b_{12r}(z)\bbT_{nr}\right)\bbD_{jr}^{-1}(z)\nonumber\\
=&~\bbL_{r}^{-1}(z)\left(\sum_{i\neq j}\beta_{ijr}(z)\bbr_{ir}\bbr_{ir}^*\bbD_{ijr}^{-1}(z)-\frac{n-1}{n}b_{12r}(z)\bbT_{nr} \bbD_{jr}^{-1}(z)\right)\nonumber\\
=&~b_{12r}(z)\bbR_{1r}(z)+\bbR_{2r}(z)+\bbR_{3r}(z),
\label{dj-1}
\end{align}
where
\begin{gather*}
\bbR_{1r}(z)=\sum_{i\neq j}\bbL_{r}^{-1}(z)(\bbr_{ir}\bbr_{ir}^*-n^{-1}\bbT_{nr})\bbD_{ijr}^{-1}(z),\\
\bbR_{2r}(z)=\sum_{i\neq j}(\beta_{ijr}(z)-b_{12r}(z))\bbL_{r}^{-1}(z)\bbr_{ir}\bbr_{ir}^*\bbD_{ijr}^{-1}(z),\\
\bbR_{3r}(z)=n^{-1}b_{12r}(z)\bbL_{r}^{-1}(z)\bbT_{nr} \sum_{i\neq j}\left(\bbD_{ijr}^{-1}(z)-\rdd{\bbD_{jr}^{-1}(z)}\right).
\end{gather*}
From this decomposition, after substituting $-\bbL_{r}^{-1}(z)$ for $\bbD_{jr}^{-1}(z)$ in the first term in \eqref{bai2.6}, there are three remaining quantities. Let's check which one (or ones) of them can be omitted.
From Lemma \ref{lem8}, \eqref{z-b}, and (4.3) of \cite{BS98}, for any $p\times p$ matrix $\bbM$, we have
\begin{align}
\rE|\tr \bbR_{2r}(z)\bbM|\leq&~ n\rE^{1/2}(|\beta_{12r}(z)-b_{12r}(z)|^2)\rE^{1/2}\bigg|\bbr_{1r}^*\bbD_{12r}^{-1}\bbM\bbL_{r}^{-1}(z)\bbr_{1r}\bigg|^2\nonumber\\
\leq&~n^{1/2} K |||\bbM|||\frac{|z|^2(1+p/(nv))}{v^5},\label{BM}
\end{align}
\rdd{where $|||\bbM|||$ denotes a nonrandom bound on the spectral norm of $\bbM$.
From Lemma \ref{lem7},}
\begin{align}
|\tr \bbR_{3r}(z)\bbM|\leq&~|||\bbM|||\frac{|z|(1+p/(nv))}{v^3}.\label{CM}
\end{align}
\rdd{Again from Lemma \ref{lem8} and \eqref{z-b}, for nonrandom $\bbM$, }
\begin{align}
\rE|\tr \bbR_{1r}(z)\bbM|\leq&~ n\rE^{1/2}|\bbr_{ir}^*\bbD_{ijr}^{-1}(z)\bbM\bbL_{r}^{-1}(z)\bbr_{ir}-n^{-1}\tr\bbT_{nr} \bbD_{ijr}^{-1}(z)\bbM\bbL_{r}^{-1}(z)|^2\nonumber\\
\leq&~n^{1/2}K||\bbM||\frac{(1+p/(nv))}{v^2}.\label{AM}
\end{align}
\rdd{Therefore, quantities containing $\bbR_{2r}(z)$ and $\bbR_{3r}(z)$ are both negligible. For the quantity involving $\bbR_{1r}(z)$, 
applying the identity $\bbD^{-1}_{js}(z)=\bbD_{ijs}^{-1}(z)-\bbD_{ijs}^{-1}(z)\bbr_{js}\bbr_{js}^*\bbD_{ijs}^{-1}(z)\beta_{ijs}(z)$,
it can be divided into three parts, that is,}
\begin{equation}\label{R123}
\tr \bbQ_r^\top\rE_j(\bbR_{1r}(z_1)) \bbQ_{r}\bbQ_s^\top \bbD_{js}^{-1}(z_2)\bbQ_{s}=R_{11}(z_1,z_2)+R_{12}(z_1,z_2)+R_{13}(z_1,z_2),
\end{equation}
where
\begin{align*}
R_{11}(z_1,z_2)&=-\sum_{i< j}\beta_{ijs}(z_2)\bbr_{ir}^*\rE_j(\bbD_{ijr}^{-1}(z_1))\bbQ_{r}\bbQ_s^\top \bbD_{ijs}^{-1}(z_2)\bbr_{is}\bbr_{is}^*\bbD_{ijs}^{-1}(z_2)\bbQ_s\bbQ_{r} ^\top\bbL_{r}^{-1}(z_1)\bbr_{ir},\\
R_{12}(z_1,z_2)&=-\tr\sum_{i<j}\bbL_{r}^{-1}(z_1)n^{-1}\bbT_{nr} \rE_j(\bbD_{ijr}^{-1}(z_1))\bbQ_r\bbQ_s^\top(\bbD_{js}^{-1}(z_2)-\bbD_{ijs}^{-1}(z_2))\bbQ_s\bbQ_r^\top,\\
R_{13}(z_1,z_2)&=\tr\sum_{i<j}\bbL_{r}^{-1}(z_1)(\bbr_{ir}\bbr_{ir}^*-n^{-1}\bbT_{nr})\rE_j(\bbD_{ijr}^{-1}(z_1))\bbQ_r\bbQ_s^\top \bbD_{ijs}^{-1}(z_2)\bbQ_s\bbQ_r^\top.
\end{align*}
From Lemma \ref{lem7} and \eqref{z-b} we get $|R_{12}(z_1,z_2)|\leq (1+p/(nv_0))/v_0^3$, and similar to \eqref{BM}, $\rE|R_{13}(z_1,z_2)|\leq n^{1/2}(1+p/(nv_0))/v_0^3$. \rdd{Thus these two parts are trivial. We then turn to dealing with $R_{11}(z_1,z_2)$.
Using Lemma \ref{lem2}, Lemma \ref{lem8}, and (4.3) of \cite{BS98} we get,} for $i<j$,
\begin{eqnarray*}
	&&\rE\bigg| \beta_{ijs}(z_2)\bbr_{ir}^*\rE_j(\bbD_{ijr}^{-1}(z_1))\bbQ_{r}\bbQ_s^\top \bbD_{ijs}^{-1}(z_2)\bbr_{is}\bbr_{is}^*\bbD_{ijs}^{-1}(z_2)\bbQ_s\bbQ_{r} ^\top\bbL_{r}^{-1}(z_1)\bbr_{ir}\\
	&&-b_{12s}(z_2)n^{-2}\tr \left(\rE_j(\bbD_{ijr}^{-1}(z_1))\bbQ_r\bbQ_s^\top \bbD_{ijs}^{-1}(z_2)\bbQ_s\bbQ_r^\top\right) \tr\left( \bbD_{ijs}^{-1}(z_2)\bbQ_s\bbQ_r^\top \bbL_{r}^{-1}(z_1)\bbQ_r\bbQ_s^\top\}\right)\bigg|\\
	&&\leq Kn^{-1/2}.
\end{eqnarray*}
\rdd{So we may simplify $R_{11}(z_1,z_2)$ by replacing $\beta_{ijs}(z_2)$ with $b_{12s}(z_2)$ and remove the random parts of $\bbr_{ir}$ and $\bbr_{is}$.
By Lemma \ref{lem7}, we have}
\begin{eqnarray*}
	&&\bigg|\tr \left(\rE_j(\bbD_{ijr}^{-1}(z_1))\bbQ_r\bbQ_s^\top \bbD_{ijs}^{-1}(z_2)\bbQ_s\bbQ_r^\top\right) \tr\left( \bbD_{ijs}^{-1}(z_2)\bbQ_s\bbQ_r^\top \bbL_{r}^{-1}(z_1)\bbQ_r\bbQ_s^\top\}\right)\\
	&&-\tr \left(\rE_j(\rdd{\bbD_{jr}^{-1}}(z_1))\bbQ_r\bbQ_s^\top \bbD_{js}^{-1}(z_2)\bbQ_s\bbQ_r^\top\right) \tr\left( \bbD_{js}^{-1}(z_2)\bbQ_s\bbQ_r^\top \bbL_{r}^{-1}(z_1)\bbQ_r\bbQ_s^\top\}\right)\bigg|\leq Kn.
\end{eqnarray*}
\rdd{It implies that we may further replace $\bbD_{ijr}^{-1}(z_1)$ and $\bbD_{ijs}^{-1}(z_2)$ in $R_{11}(z_1,z_2)$ with $\bbD_{jr}^{-1}(z_1)$ and $\bbD_{js}^{-1}(z_2)$, respectively, which yields}
\begin{eqnarray}
&&\rE\bigg|R_{11}(z_1,z_2)\nonumber\\
&&+\frac{j-1}{n^2}b_{12s}(z_2)\tr\left(\bbQ_r^\top \rE_j(\bbD_{jr}^{-1}(z_1))\bbQ_r\bbQ_s^\top \bbD_{js}^{-1}(z_2)\bbQ_s\right)\tr\left( \bbQ_s^\top\bbD_{js}^{-1}(z_2)\bbQ_s\bbQ_r^\top \bbL_{r}^{-1}(z_1)\bbQ_r\right)\bigg|\nonumber\\
&&\leq Kn^{1/2}.\label{az12}
\end{eqnarray}
\rdd{Integrating the results in \eqref{dj-1}-\eqref{az12}, we obtain}
\begin{eqnarray}
	&&\tr\left( \bbQ_r^\top\rE_j(\bbD_{jr}^{-1}(z_1))\bbQ_r\bbQ_s^\top \bbD_{js}^{-1}(z_2)\bbQ_s \right)\non
	&&\times\left(1+\frac{j-1}{n^2}b_{12r}(z_1)b_{12s}(z_2)\tr\left( \bbQ_s^\top\bbD_{js}^{-1}(z_2)\bbQ_s\bbQ_r^\top \bbL_{r}^{-1}(z_1)\bbQ_r\right)\right)\non
	&=&-\tr \bbQ_r^\top\bbL_{r}^{-1}(z_1)\bbQ_r\bbQ_s^\top \bbD_{js}^{-1}(z_2)\bbQ_s +R_{14}(z_1,z_2),\label{f2}
\end{eqnarray}
where $\rE|R_{14}(z_1,z_2)|\leq Kn^{1/2}$. Furthermore, from this and \eqref{dj-1}-\eqref{AM}, \rdd{we may substitute for the second and third $\bbD_{js}^{-1}(z_2)$ in \eqref{f2} with $-\bbL_s^{-1}(z_2)$ and then get
} 
\begin{eqnarray}
&&\tr\left( \bbQ_r^\top\rE_j(\bbD_{jr}^{-1}(z_1))\bbQ_r\bbQ_s^\top \bbD_{js}^{-1}(z_2)\bbQ_s \right)\non
&&\times\left(1-\frac{j-1}{n^2}b_{12r}(z_1)b_{12s}(z_2)\tr\left( \bbQ_s^\top\bbL_{s}^{-1}(z_2)\bbQ_s\bbQ_r^\top \bbL_{r}^{-1}(z_1)\bbQ_r\right)\right)\nonumber\\
&=&\tr \bbQ_r^\top\bbL_{r}^{-1}(z_1)\bbQ_r\bbQ_s^\top \bbL_{s}^{-1}(z_2)\bbQ_s +R_{15}(z_1,z_2),\label{R15}
\end{eqnarray}
where $\rE|R_{15}(z_1,z_2)|\leq Kn^{1/2}$.

From Lemma \ref{lem7} and (4.3) of \cite{BS98}, we have
$$
|b_{12r}(z)-b_{nr}(z)|\leq Kn^{-1}\quad \text{and}\quad |b_{nr}(z)-\rE\beta_{1r}(z)|\leq Kn^{-1/2},
$$
respectively. By (2.2) of \cite{Silv95} and discussions in Section 5 of \cite{BS98}, we have
$$
\rE\beta_{1r}(z)=-z\rE\underline{m}_{nr}(z)\quad\text{and}\quad |\rE\um_{nr}(z)-\underline{m}_{nr}^0(z)|\leq Kn^{-1},
$$
respectively. Therefore, we get
\begin{equation}\label{b-zm}
|b_{12r}(z)+z\underline m_{nr}^0(z)|\leq Kn^{-1/2}.
\end{equation}
\rdd{Combining this and \eqref{R15}}, it follows that
\begin{eqnarray}
&&\tr\left( \bbQ_r^\top \rE_j(\bbD_{jr}^{-1}(z_1))\bbQ_r\bbQ_s^\top \bbD_{js}^{-1}(z_2)\bbQ_s\right)\left(1-\frac{j-1}{n^2}\um_{nr}^0(z_1)\um_{ns}^0(z_2)\right.\nonumber\\
&&\times\tr\left( \bbQ_s^\top(I+\um_{ns}^0(z_2)\bbT_{ns})^{-1}\bbQ_s\bbQ_r^\top  (I+\um_{nr}^0(z_1)\bbT_{nr})^{-1}\bbQ_r\right)\bigg)\nonumber\\
&=&\tr\left( \bbQ_r^\top \rE_j(\bbD_{jr}^{-1}(z_1))\bbQ_r\bbQ_s^\top \bbD_{js}^{-1}(z_2)\bbQ_s\right)\left(1-\frac{j-1}{n}c_n\um_{nr}^0(z_1)\um_{ns}^0(z_2)\right.\nonumber\\
&&\times\int\int\frac{t_1t_2dH_{prs}(t_1,t_2)}{(1+t_1\um_{nr}^0(z_1))(1+t_2\um_{ns}^0(z_2))}\bigg)\nonumber\\
&=&\frac{nc_n}{z_1z_2}\int\int\frac{t_1t_2dH_{prs}(t_1,t_2)}{(1+t_1\um_{nr}^0(z_1))(1+t_2\um_{ns}^0(z_2))} +R_{16}(z_1,z_2),\label{var-2}
\end{eqnarray}
where $\rE|R_{16}(z_1,z_2)|\leq Kn^{1/2}$. Similar to the arguments in \cite{BS04} (page 577),  \rdd{using \eqref{var-2}} and letting
$$
a_n(z_1,z_2)=c_n\um_{nr}^0(z_1)\um_{ns}^0(z_2)\int\int\frac{t_1t_2dH_{prs}(t_1,t_2)}{(1+t_1\um_{nr}^0(z_1))(1+t_2\um_{ns}^0(z_2))},
$$
we get
\begin{eqnarray}\label{part1}
\Gamma_1=\frac{1}{n}\sum_{j=1}^n\frac{a_n(z_1,z_2)}{1-((j-1)/n)a_n(z_1,z_2)}+o_p(1)\xrightarrow{i.p.}-\log(1-a(z_1,z_2)),
\end{eqnarray}
where
\begin{eqnarray*}
	a(z_1,z_2)=c\um_r(z_1)\um_s(z_2)\int\int\frac{t_1t_2dH_{rs}(t_1,t_2)}{(1+t_r\um_r(z_1))(1+t_s\um_s(z_2))}.
\end{eqnarray*}
Similar to the derivation of $\Gamma_1$, one can easily show that
\begin{eqnarray}\label{part2}
\alpha_x\Gamma_2=\frac{1}{n}\sum_{j=1}^n\frac{\alpha_xa_n(z_1,z_2)}{1-\alpha_x((j-1)/n)a_n(z_1,z_2)}+o_p(1)\xrightarrow{i.p.}-\log(1-\alpha_xa(z_1,z_2)).
\end{eqnarray}
\rdd{Considering the third term of \eqref{bai2.6}, $\beta_x\Gamma_3$ with $\beta_x\neq 0$,} from Assumption (e), the matrix \rdd{$\bbQ_r^\top\bbQ_r$ is diagonal, so is $\bbQ_r^\top\bbL_r^{-1}\bbQ_r$.} \rdd{Using \eqref{dj-1}-\eqref{AM} and \eqref{b-zm}, we have
\begin{align}
\beta_x \Gamma_3=&~\beta_xb_{nr}(z_1)b_{ns}(z_2)\frac{1}{n^2}\sum_{j=1}^{n}\sum\limits_{i=1}^{k}{\bf e}_i^\top\bbQ_r^\top\rE_j(\bbD_{jr}^{-1}(z_1))\bbQ_r{\bf e}_i{\bf e}_i^\top\bbQ_s^\top\rE_j(\bbD_{js}^{-1}(z_2))\bbQ_s{\bf e}_i\nonumber\\
=&~\beta_xb_{nr}(z_1)b_{ns}(z_2)\frac{1}{n^2}\sum_{j=1}^{n}\sum\limits_{i=1}^{k}{\bf e}_i^\top\bbQ_r^\top\bbL_r^{-1}(z_1)\bbQ_{r}{\bf e}_i{\bf e}_i^\top\bbQ_s^\top\bbL_s^{-1}(z_2)\bbQ_s{\bf e}_i+o_p(1)\nonumber\\
=&~\beta_xb_{nr}(z_1)b_{ns}(z_2)\frac{1}{n^2}\sum_{j=1}^{n}\tr\left(\bbQ_r^\top\bbL_r^{-1}(z_1)\bbQ_{r}\bbQ_s^\top\bbL_s^{-1}(z_2)\bbQ_s\right)+o_p(1)\nonumber\\
\xrightarrow{i.p.}& ~\beta_xa(z_1,z_2).\label{part3}
\end{align}}

Collecting the results in \eqref{limcov}, \eqref{part1}-\eqref{part3},
we finally get the covariance function in the Lemma and the proof of
this Lemma is completed.

\subsection{Proof of Theorem \ref{MainThmMult}}

\rvv{
First we show that it is 
it is enough to establish the following claim:   under the (classical)  Mar\v{c}enko-Pastur regime, i.e., 
$ n\to \infty$, $p=p_n\to \infty$ such that  $p/n\to c>0$, it holds
that 
\begin{equation}
  \label{eq:claim}
  \phi_{q,n} - \frac{q}{2} 
          \xrightarrow{d} \N\lb 0, s(c) \rb,
\end{equation}
where recall that $s(u)=  q+u(\nu_4-1)(q^2+\frac{q}{2})$.
Here we use $\phi_{q,n}$ for $\phi_q$ to signify the dependence in
$n$. 
So assume  this claim is true.  Under the relaxed Mar\v{c}enko-Pastur
regime, the sequence $\{p_n/n\}$ is bounded below and above. 
For any  subsequence  $(p_{n_k}, n_k)_k$ of $(p_n,n)$,
we can extract a further subsequence  
$(p_{n_{k_\ell}}, n_{k_\ell})_\ell$ such that the ratios 
$  p_{n_{k_\ell}} /  n_{k_\ell}   $  converge to $\alpha >0$  when
$\ell\to\infty$.
On this subsequence, by Claim~\eqref{eq:claim},  
\[  
\phi_{q, n_{k_\ell}} - \frac{q}{2} 
\xrightarrow{d} \N\lb
0,  s(\alpha) \rb, \quad \ell\to\infty.
 \]
By continuity of the function $u\to s(u)$, we have 
\[   
s(  p_{n_{k_\ell}} / n_{k_\ell} )^{-1/2}  \lb \phi_{q, n_{k_\ell}} - \frac{q}{2}\rb 
\xrightarrow{d} \N\lb 0, 1 \rb, \quad \ell\to\infty.
\] 
As this limit is independent of the subsequence and it holds for all
such subsequences, the same limit holds for the whole sequence, that
is, 
\[   
s(  p_{n} / n )^{-1/2}  \lb \phi_{q, n} - \frac{q}{2}\rb 
\xrightarrow{d} \N\lb 0, 1 \rb, \quad \ell\to\infty.
\] 
The required asymptotic normality is thus established. 

The remaining of the section is devoted to a proof of Claim
\eqref{eq:claim} assuming $p/n\to c>0$.
}
Define the banded Toeplitz matrix
\[\bbC_{n,\tau}=
  \left(
                                            \begin{array}{ccccc}
                                              0 & \cdots & \frac{1}{2} & \cdots & 0 \\
                                              \vdots & \ddots & 0 & \frac12 & \vdots \\
                                              \frac12 & 0 & \ddots & 0 & \frac12 \\
                                              \vdots & \frac12 & 0 & \ddots & \vdots \\
                                              0 & \cdots & \frac12 & \cdots & 0 \\
                                            \end{array}
                                          \right)_{n\times n}
\]
and
\[
  \widehat{\bbN}_{\tau}=\frac {1}{2p}\sum_{t=1+\tau}^{n}\lb\bbx_t\bbx_{t-\tau}^*+\bbx_{t-\tau}\bbx_t^*\rb=\frac1p \bbX_n\bbC_{n,\tau}\bbX_{n}^*.
\]

Define the associated  Fourier series $f(\lambda)$ of the banded
Toeplitz matrix $\bbC_{n,\tau}$ as
\begin{align*}
f(\lambda)&=\lim_{n\rightarrow \infty}\sum_{k=-\tau}^\tau t_k e^{i k \lambda}= \frac 12\lb e^{i\tau\lambda}+e^{-i\tau\lambda}\rb=\cos(\tau\lambda),
\end{align*}
where $t_k$ is entry on the $k$-diagonal of $\bbC_{n,\tau}$.
\par
\vspace{0.5cm}
According to the fundamental eigenvalue
distribution theorem of Szeg\"{o} for Toeplitz forms, see Section 5.2 in \citet{Grenander58} and Theorem 4.1 in \citet{Gray06}, we can infer that
\begin{lemma}\label{lemma2}
  Suppose $\{l_t,~t=1,\cdots,n\}$ are eigenvalues of $\bbC_{n,\tau}$ with Fourier series $f(\lambda)$, then
  \begin{itemize}
  \item[(1)]
  For any positive integer $s$,
  \[\lim_{n\rightarrow \infty}\frac 1n\sum_{t=1}^n l_t^s=\frac{1}{2\pi}\int_0^{2\pi}f(\lambda)^s\d \lambda=\frac{1}{2\pi}\int_0^{2\pi}\lb\cos(\tau\lambda)\rb^s \d\lambda, \]
  \item[(2)]
  For any continuous function on support of $\{l_t,~t=1,\cdots,n\}$,
  \begin{equation*}
    \lim_{n\rightarrow \infty}\frac 1n\sum_{t=1}^n F(l_t)=\frac{1}{2\pi}\int_0^{2\pi} F\lb\cos(\tau\lambda)\rb\d \lambda,
  \end{equation*}
 \item[(3)] Sequence $\left\{l_t,~t=1,\cdots,n\right\}$ and
 \begin{equation*}
   \left\{\cos\lb \frac{2\pi \tau t}{n}\rb,~t=1,\cdots,n\right\}
 \end{equation*} are asymptotically equally distributed.
  \end{itemize}
\end{lemma}

The limiting spectral distribution of $\bbC_{n,\tau}$ is also derived in Lemma 3.1 of \citet{Bai15}, this useful lemma is stated as follows:

\begin{lemma}
 As $T\rightarrow \infty$, \rdd{the ESD} of $\bbC_{n,\tau}$ tends to $H$, which is
an Arcsine distribution with density function
\[H'(t)=\cfrac{1}{\pi \sqrt{1-t^2}},~t\in(-1,1).\]
\end{lemma}

\rddd{Recall for the  permutation matrices
$\bbD_1$ and $\bbD_{\tau}=\bbD_1^{\tau}$ defined in Introduction, it holds that}
\[
  \bbD_1\bbD_1^\top=\bbD_1^\top\bbD_1=\bbD_\tau
  \bbD_\tau^\top=\bbD_\tau^\top\bbD_\tau={\bf I}_n.
\]
Meanwhile, from the properties of Chebyshev polynomials, we can derive the following lemma.
\begin{lemma}\label{lem3}
  \begin{itemize}
    \item[(1)] $\frac{1}{2}\lb \bbD_1+\bbD_1^\top\rb$ has eigenvalue $$\left\{\cos\lb\frac{2\pi t}{n}\rb,~t=1,\cdots,n \right\},$$
    \item[(2)] $\frac{1}{2}\lb \bbD_{\tau}+\bbD_{\tau}^\top\rb$ has eigenvalue
    \[\left\{T_\tau\lb \cos\lb\frac{2\pi t}{n}\rb \rb,~t=1,\cdots,n\right\},\]
    where $T_{\tau}(\cdot)$ stands for the Chebyshev polynomial of order $\tau$.
    \item[(3)]$\frac{1}{2}\lb \bbD_{\tau}+\bbD_{\tau}^\top\rb$ shares the same asymptotic spectral distribution with $\bbC_{n,\tau}$ as $n\rightarrow\infty$.
  \end{itemize}
\end{lemma}

Since
\begin{align*}
 \widetilde{\bbN}_{\tau}&=\frac {1}{2p}\sum_{t=1}^{n}\lb\bbx_t\bbx_{t-\tau}^*+\bbx_{t-\tau}\bbx_t^*\rb\\
 &=\frac{1}{2p}\lb\bbx_1,\cdots,\bbx_n\rb\lb\sum_{\tau=1}^q\lb \bbD_{\tau}+\bbD_{\tau}^\top\rb\rb \lb\bbx_1,\cdots,\bbx_n\rb^*\\
  &=\frac{1}{2p} \bbX_n\lb\sum_{\tau=1}^q\lb \bbD_{\tau}+\bbD_{\tau}^\top\rb \rb \bbX_n^*,
\end{align*}
here for $t\leq 0$, $\bbx_t=\bbx_{n+t}$, by \rdd{Lemmas} \ref{lemma2} and \ref{lem3}, it doesn't take too much effort to see that $\widetilde{\bbN}_{\tau}$ and $\widehat{\bbN}_\tau$ share the same limiting spectral distribution.

Consider the Stieltjes transform $\underline{m}_\tau(z)$ of the limiting spectral distribution of $\widetilde{\bbN}_\tau$, by implementing the Silverstein equation \eqref{eq:Silv}, we can infer  that $\um_\tau(z)$ satisfies
\begin{align*}
    z&=-\frac1{\um_\tau(z)} +\frac{1}{c}\int\frac{t}{1+t\um_\tau(z)}dH_\tau(t)\\
    &=\frac{1}{\um_{\tau}(z)}\lb-1+\frac{1}{c}-\frac{1}{c
    	\sqrt{1-\um_{\tau}^2(z)}}\rb,    
\end{align*}
where \rdd{$p/n\rightarrow c>0$} as $n\rightarrow \infty$, which coincides with the results in \citet{Bai15} and \citet{Li16}.

Note that our test statistic
$$\mathcal{L}_q=\sum_{\tau=1}^q\widetilde{L}_\tau=\sum_{\tau=1}^q\mathrm{Tr}(\widetilde{\bbM}_{\tau}\widetilde{\bbM}_{\tau}^*)=\lb\frac{p}{n}\rb^2\sum_{\tau=1}^q\mathrm{Tr}(\widetilde{\bbN}_{\tau}\widetilde{\bbN}_{\tau}^*),$$
where $  \widetilde{\bbM}_{\tau}= \frac12
  \left(\widehat{\bbSigma}_{\tau} +\widehat{\bbSigma}_{\tau}^* \right)
  =\frac
  {1}{2n}\sum_{t=1}^{n}\lb\bbx_t\bbx_{t-\tau}^*+\bbx_{t-\tau}\bbx_t^*\rb$,
thus the asymptotic properties of $\widetilde{\bbM}_\tau$ can be inferred from those of $\widetilde{\bbN}_{\tau}$ since $\widetilde{\bbM}_\tau=\frac{p}{n}\widetilde{\bbN}_\tau$.

Actually in \citet{Li16}, the asymptotic behavior of the single-lag-$\tau$ test statistic $\widetilde{L}_\tau$  has already been thoroughly explored and characterized. Theorem 2.1 in \citet{Li16} is stated as follows:

\begin{lemma}\label{MainThmSig}
  Let $\tau\geq 1$ be a fixed integer, and assume that
  \begin{enumerate}
  \item
    $\{z_{it},~i=1,\cdots,p,~t=1,\cdots,n\}$ are all independently distributed satisfying
    $\rE z_{it}=0,~\rE z_{it}^2=1,~\rE z_{it}^4=\nu_4<\infty$;
  \item (Mar\v{c}enko-Pastur  regime).
    The dimension $p$ and the sample size $n$ grow  to infinity in
    a related way such that     $c_n:=p/n\rightarrow c>0$.
  \end{enumerate}
  Then in the simplest setting when $\bbx_t=\z_t$, the limiting distribution of the test statistic $\widetilde{L}_\tau$ is
  \begin{equation*}
    \frac{n}{p}\widetilde{L}_\tau-\frac p2 \xrightarrow{d} \N\lb \frac{1}{2},~1+\frac{3(\nu_4-1)}{2}c\rb.
  \end{equation*}
\end{lemma}

Now consider the multi-lag-$q$ test statistic $\mathcal{L}_q=\sum_{\tau=1}^q\widetilde{L}_\tau$, combining with Lemma \ref{MainThmSig}, all we need is the joint distribution of any two different single-lag test statistic, i.e.
$\lb\widetilde{L}_r, \widetilde{L}_s\rb,~1\leq r\neq s\leq q$.

For a given integer $q>0$, $1\leq r\neq s\leq q$, let $f_r(x)=f_s(x)=x^2,$
\begin{align*}
    \bbB_{nr}&=\frac{1}{2p}\lb \bbD_{r}+\bbD_{r}^\top\rb \bbX_n^* \bbX_n,~B_{ns}=\frac{1}{2p}\lb \bbD_{s}+\bbD_{s}^\top\rb \bbX_n^*\bbX_n,\\
  \underline{\bbB}_{nr}&=\widetilde{\bbN}_r=\frac{1}{2p} \bbX_n\lb \bbD_{r}+\bbD_{r}^\top\rb \bbX_n^*,
  ~\underline{B}_{ns}=\widetilde{\bbN}_s=\frac{1}{2p} \bbX_n\lb \bbD_{s}+\bbD_{s}^\top\rb \bbX_n^*,
\end{align*}
then both the \rdd{LSDs} of $\underline{\bbB}_{nr}$ and
$\underline{\bbB}_{ns}$,
i.e. $\underline{F}^{c,H_r},~\underline{F}^{c,H_s}$ have Stieltjes
transform $\um_r(z)$ and $\um_s(z)$ satisfying the equation
\[z=\frac{1}{\um(z)}\lb-1+\frac{1}{c}-\frac{1}{c\sqrt{1-\um^2(z)}}\rb.
\]
Meanwhile,
$$\mathrm{Tr}(\widetilde{\bbN}_{r}\widetilde{\bbN}_{r}^*)=\int f_r(x)\d F_{nr}(x),~\mathrm{Tr}(\widetilde{\bbN}_{s}\widetilde{\bbN}_{s}^*)=\int f_s(x)\d F_{ns}(x).$$
where $F_{nr}$ and $F_{ns}$ are the \rdd{ESDs} of $\bbB_{nr}$ and
$\bbB_{ns}$. Thus, by directly implementing our joint CLT for linear
spectral statistics of the sample covariance matrices, i.e. Theorem
\ref{thm2}, we can derive the joint distribution of
$\lb\widetilde{L}_1,\cdots, \widetilde{L}_q\rb$ or any pair of
$\lb\widetilde{L}_r, \widetilde{L}_s\rb,~1\leq r\neq s\leq q$.
\rddd{Precisely, the covariance function in Theorem \ref{thm2} for the
present case can be calculated to be 
\begin{equation}
  \label{eq:appcov}
  \renewcommand{\arraystretch}{1.5}
  \Cov\left(X_{f_{r}},X_{f_{s}}\right)=\left\{
    \begin{array}{ll}
      \frac{1+\frac{3}{2}c(\beta_x+2)}{c^2}, &~\mbox{ if } r=s,\\
      \frac{\beta_x+2}{c}, &~\mbox{ if } r\neq s.
    \end{array}
  \right.
\end{equation}
The details of this lengthy derivation are postponed to Appendix~\ref{sec:appcov}.}
Combining with Lemma \ref{MainThmSig}, for any given integer $q\leq 1$, it can be inferred that, under the same assumptions in Lemma \ref{MainThmSig}, the joint limiting distribution of  $\lb \widetilde{L}_1, \cdots, ~\widetilde{L}_q\rb$ is
\[
    \frac{n}{p}\lb
    \begin{array}{c}
      \widetilde{L}_1\\
      \vdots\\
      \widetilde{L}_q
    \end{array}
    \rb-\frac p2 {\bf 1}_q \xrightarrow{d} \N\lb \frac{1}{2}{\bf 1}_q,~
        \begin{pmatrix}
      1+\frac{3(\nu_4-1)}{2}c&\cdots&c(\nu_4-1)\\
      \vdots&\ddots&\vdots\\
      c(\nu_4-1)&\cdots&1+\frac{3(\nu_4-1)}{2}c
    \end{pmatrix}
    \rb,
\]
where ${\bf 1}_q=\lb 1,\cdots,1\rb^\top$ is a $q$ dimensional vector with $q$ ones.

Recall that
$$\mathcal{L}_q=\sum_{\tau=1}^q\widetilde{L}_\tau=\lb\frac{p}{n}\rb^2\sum_{\tau=1}^q\mathrm{Tr}(\widetilde{\bbN}_{\tau}\widetilde{\bbN}_{\tau}^*),$$
then by \rdd{the} Delta method, we can derive the limiting distribution of our test statistic $\mathcal{L}_q$, i.e.,
 \[\frac{n}{p}\mathcal{L}_q-\frac{pq}{2} \xrightarrow{d} \N\lb \frac{q}{2},~q+c(\nu_4-1)(q^2+\frac{q}{2})\rb.\]
\rvv{Claim \eqref{eq:claim} is thus established.}

\section{Discussions}

In this paper we have introduced, for the first time in the literature
on eigenvalues of large sample covariance matrices, a joint central
limit theorem  involving several population covariance matrices. This
theorem is believed to provide wide applications to current
problems in high-dimensional statistics, especially for testing on
structures of population covariances. As a show-case, we treated the
problem of testing for a high-dimensional white noise in time series
modelling.
\rddd{The derived new test shows very promising performance compared
  to existing competitors}
For future study, it would be worth investigating other
significant applications of this CLT.

\medskip 
\noindent
{\bf Acknowledgments.} \quad
\rddd{We thank a Referee for the suggestion of the comparison to a
permutation test made in Section 3.1.}
{Weiming Li's research is partially supported by National Natural Science Foundation of China, No. 11401037, MOE Project of Humanities 
  and Social Sciences, No. 17YJC790057, and Program of IRTSHUFE.}
Jianfeng Yao thanks support from HKSAR GRF Grant 17305814.

\appendix
\section{Mathematical Tools}



\begin{lemma}\label{lem2}
  \rdd{Let $\bbX=(X_1,\ldots,X_n)^\top$ be a (complex) random vector with independent and standardized entries having finite fourth moment and $\bbC=(c_{ij})$ be $n\times n$ (complex) matrix. We have}
  $$\rE|\bbX^*\bbC\bbX-\rtr\bbC|^4\leq K\left(\left[
      \rtr  (\bbC\bbC^*)\right]^{2}
    +\siln \rE|X_{ii}^8||c_{ii}|^4\right). $$
\end{lemma}
The proof of the lemma \rdd{follows easily} by simple calculus and thus omitted.



\begin{lemma} \label{lem7} (Lemma 2.6 of \citet{SB95}).
  Let $z\in {\mathbb C}^+$ with $v=\Im\,z$, $\bbA$ and $\bbB$ being $n\times n$ with $\bbB$
  Hermitian, and $\bbr\in{\mathbb C} ^n$.  Then
  $$\bigl|\rtr\bigl((\bbB-z\bbI)^{-1}-(\bbB+\bbr\bbr^*-z\bbI)^{-1}\bigr)\bbA\bigr|=
  \left|\frac{\bbr^*( \bbB-z\bbI)^{-1}\bbA( \bbB-z\bbI)^{-1}\bbr}
    {1+\bbr^*(\bbB-z\bbI)^{-1}\bbr}\right|\leq\frac{\|\bbA\|}v.$$
\end{lemma}

\begin{lemma}\label{lem8}[Formula 2.3 of \citet{BS04}]
For any nonrandom $p\times p$ matrices $\bbC_k$, $k=1,\ldots,q_1$ and $\tilde \bbC_\ell$, $\ell=1,\ldots q_2$.
\begin{eqnarray}
&&\bigg|\rE\left(\prod_{k=1}^{q_1}\bbr_{1r}^*\bbC_k\bbr_{1r}\prod_{\ell=1}^{q_2}(\bbr_{1r}^*\tilde \bbC_\ell\bbr_{1r}-n^{-1}\tr \bbT_{nr}\tilde \bbC_\ell)\right)\bigg|\nonumber\\
&\leq& Kn^{-(1\wedge q_2)}\delta_n^{(q_2-2)\vee0}\prod_{k=1}^{q_1}||\bbC_k||\prod_{\ell=1}^{q_2}||\tilde \bbC_\ell ||,\quad q_1, q_2\geq0,\label{q12-moment}
\end{eqnarray}
where $K$ is a positive constant depending on $q_1$ and $q_2$.
\end{lemma}

{\color{red}
\section{Derivation of the covariance \protect\eqref{eq:appcov}}
\label{sec:appcov}
}

Applying Theorem \ref{thm2} to the functions $f_r(x)=f_s(x)$ where
$1\le r\ne s\le q$, 
the corresponding  covariance function is
\begin{eqnarray}
    {\rm Cov}(X_{f_{r}}, X_{f_{s}})
    &=&\frac{1}{4\pi^2}\oint\limits_{{\cal
        C}_1}\oint\limits_{{\cal
        C}_2} f_{r}(z_1)f_{s}(z_2)\frac{\partial^2g(z_1,z_2)}{\partial z_1\partial z_2}dz_1dz_2,
  \end{eqnarray}
  where $g(z_1,z_2)=\log(1-a(z_1,z_2))+\log(1-\alpha_xa(z_2,z_2))-\beta_xa(z_1,z_2)$ with
  \begin{eqnarray*}
  a(z_1,z_2)=\iint\frac{c\um_r(z_1)\um_s(z_2)t_1t_2}{(1+t_1\um_r(z_1))(1+t_2\um_s(z_2))}dH_{rs}(t_1,t_2).
 \end{eqnarray*}
Mapping into our case, we have
\[p\leftrightarrow n, ~n\leftrightarrow p,~c\leftrightarrow \frac 1c,~\alpha_x=1,~\beta_x=\nu_4-3,~f_{r}(x)=f_{s}(x)=x^2,\]
\begin{align*}
  a(z_1,z_2)&=\frac{1}{c}\iint\frac{\um_r(z_1)\um_s(z_2)t_1t_2}{(1+t_1\um_r(z_1))(1+t_2\um_s(z_2))}dH_{rs}(t_1,t_2)\\
&=\frac{1}{c}\int_{-1}^1 \cfrac{\um_1\um_2 T_r(t)T_s(t)}{(1+T_r(t)\um_1)(1+T_s(t)\um_2)}\d H(t),
\end{align*}
where $\um_1\triangleq\um_r(z_1),~\um_2\triangleq\um_s(z_2)$, $T_r(t),~ T_s(t)$ are Chebyshev polynomial of order $r$ and $s$ respectively. Furthermore, we have

\begin{align*}
  \frac{\partial a\left(z_{1},z_{2}\right)}{\partial z_{1}}
   & =  \frac{1}{c}\int_{-1}^1\frac{T_{s}\left(t\right)\underline{m}_{2}}{1+T_{s}\left(t\right)\underline{m}_{2}}\cdot\frac{T_{r}(t)}{\left(1+T_{r}\left(t\right)\underline{m}_{1}\right)^{2}}\cdot\frac{\partial\underline{m}_{1}}{\partial z_{1}}{\rm d}H(t),
 \end{align*}
 \begin{gather*}
\frac{\partial a\left(z_{1},z_{2}\right)}{\partial z_{2}}  = \frac{1}{c}\int_{-1}^{1}\frac{T_{r}\left(t\right)\underline{m}_{1}}{1+T_{r}\left(t\right)\underline{m}_{1}}\cdot\frac{T_{s}(t)}{\left(1+T_{s}\left(t\right)\underline{m}_{2}\right)^{2}}\cdot\frac{\partial\underline{m}_{2}}{\partial z_{2}}{\rm d}H(t),\\
\frac{\partial^{2}a\left(z_{1},z_{2}\right)}{\partial z_{1}\partial z_{2}}  = \frac{1}{c}\int_{-1}^{1}\frac{T_{r}(t)}{\left(1+T_{r}\left(t\right)\underline{m}_{1}\right)^{2}}\cdot\frac{T_{s}(t)}{\left(1+T_{s}\left(t\right)\underline{m}_{2}\right)^{2}}\cdot\frac{\partial\underline{m}_{1}}{\partial z_{1}}\cdot\frac{\partial\underline{m}_{2}}{\partial z_{2}}{\rm \mathrm{d}}H(t).
\end{gather*}

Since $g(z_{1},z_{2})=2\log\left(1-a(z_{1},z_{2})\right)-\beta_{x}a(z_{1},z_{2}),$
\[
\frac{\partial^{2}\log\left(1-a(z_{1},z_{2})\right)}{\partial z_{1}\partial z_{2}}=-\frac{\frac{\partial^{2}a}{\partial z_{1}\partial z_{2}}}{1-a}-\cfrac{\frac{\partial a}{\partial z_{1}}\cdot\frac{\partial a}{\partial z_{2}}}{\left(1-a\right)^{2}},\]
Thus
\begin{align*}
\Cov\left(X_{f_{r}},X_{f_{s}}\right) & = \frac{1}{4\pi^{2}}\oint_{\mathcal{C}_{1}}\oint_{\mathcal{C}_{2}}z_{1}^{2}z_{2}^{2}\frac{\partial^{2}g\left(z_{1},z_{2}\right)}{\partial z_{1}\partial z_{2}}\mathrm{d}z_{1}\mathrm{d}z_{2}\\
 & =  -\frac{1}{2\pi^{2}}\oint_{\mathcal{C}_{1}}\oint_{\mathcal{C}_{2}}z_{1}^{2}z_{2}^{2}\frac{\frac{\partial^{2}a}{\partial z_{1}\partial z_{2}}}{1-a}\mathrm{d}z_{1}\mathrm{d}z_{2}-\frac{1}{2\pi^{2}}\oint_{\mathcal{C}_{1}}\oint_{\mathcal{C}_{2}}z_{1}^{2}z_{2}^{2}\cfrac{\frac{\partial a}{\partial z_{1}}\cdot\frac{\partial a}{\partial z_{2}}}{\left(1-a\right)^{2}}\mathrm{d}z_{1}\mathrm{d}z_{2}\\
 & ~ -\frac{\beta_{x}}{4\pi^{2}}\oint_{\mathcal{C}_{1}}\oint_{\mathcal{C}_{2}}z_{1}^{2}z_{2}^{2}\frac{\partial^{2}a}{\partial z_{1}\partial z_{2}}\mathrm{d}z_{1}\mathrm{d}z_{2}\triangleq M_{1}+M_{2}+M_{3}.
 \end{align*}

Consider $M_{1}$ first, by Cauchy's residue theorem, we have

\begin{align*}
~&~\frac{1}{2\pi i}\oint_{\mathcal{C}_{1}}\frac{z_{1}^{2}}{1-a}\cdot\frac{\partial^{2}a}{\partial z_{1}\partial z_{2}}\mathrm{d}z_{1} \\
 = & ~\frac{1}{2\pi i}\oint_{\mathcal{C}_{1}}\frac{z_{1}^{2}}{1-a}\cdot\frac{1}{c}\int_{-1}^{1}\frac{T_{r}(t)}{\left(1+T_{r}\left(t\right)\underline{m}_{1}\right)^{2}}\cdot\frac{T_{s}(t)}{\left(1+T_{s}\left(t\right)\underline{m}_{2}\right)^{2}}{\rm \mathrm{d}}H(t)\cdot\frac{\partial\underline{m}_{1}}{\partial z_{1}}\cdot\frac{\partial\underline{m}_{2}}{\partial z_{2}}\mathrm{d}z_{1}\\
  = & -\frac{\partial\underline{m}_{2}}{\partial z_{2}}\cdot\frac{1}{c}\int_{-1}^{1}\left[\frac{T_{r}(t)T_{s}(t)}{\left(1+T_{s}\left(t\right)\underline{m}_{2}\right)^{2}}\cdot\frac{1}{2\pi i}\oint_{\mathcal{C}_{1}}\frac{z_{1}^{2}}{\left(1-a\right)\left(1+T_{r}\left(t\right)\underline{m}_{1}\right)^{2}}\mathrm{d}\underline{m}_{1}\right]\mathrm{d}H(t)\\
  = & -\frac{\partial\underline{m}_{2}}{\partial z_{2}}\cdot\frac{1}{c}\int_{-1}^{1}\left[\frac{T_{r}(t)T_{s}(t)}{\left(1+T_{s}\left(t\right)\underline{m}_{2}\right)^{2}}\cdot\frac{1}{2\pi i}\oint_{\mathcal{C}_{1}}\frac{\left(-1+\frac{1}{c}-\frac{1}{c\sqrt{1-\underline{m}_{1}^{2}}}\right)^{2}}{\underline{m}_{1}^{2}\left(1-a\right)\left(1+T_{r}\left(t\right)\underline{m}_{1}\right)^{2}}\mathrm{d}\underline{m}_{1}\right]\mathrm{d}H(t)\\
 = & -\frac{\partial\underline{m}_{2}}{\partial z_{2}}\cdot\frac{1}{c}\int_{-1}^{1}\left[\frac{T_{r}(t)T_{s}(t)}{\left(1+T_{s}\left(t\right)\underline{m}_{2}\right)^{2}}\cdot\left.\left[\frac{\left(-1+\frac{1}{c}-\frac{1}{c\sqrt{1-\underline{m}_{1}^{2}}}\right)^{2}}{\left(1-a\right)\left(1+T_{r}\left(t\right)\underline{m}_{1}\right)^{2}}\right]^{(1)}\right\rvert_{\um_1=0} \right]\mathrm{d}H(t).
 \end{align*}

\noindent
Note that
\[
\frac{\partial a}{\partial\underline{m}_{1}}=\frac{1}{c}\int_{-1}^{1}\frac{T_{s}(t)\underline{m}_{2}}{1+T_{s}\left(t\right)\underline{m}_{2}}\cdot\frac{T_{r}(t)}{\left(1+T_{r}\left(t\right)\underline{m}_{1}\right)^{2}}\mathrm{d}H(t),\]

\noindent
then
\[
\left.\left[\frac{\left(-1+\frac{1}{c}-\frac{1}{c\sqrt{1-\underline{m}_{1}^{2}}}\right)^{2}}{\left(1-a\right)\left(1+T_{r}\left(t\right)\underline{m}_{1}\right)^{2}}\right]^{(1)}\right\rvert_{\um_1=0} =-2T_{r}(t)+\frac{1}{c}\int_{-1}^{1}\frac{T_{s}(u)T_{r}(u)\underline{m}_{2}}{1+T_{s}\left(u\right)\underline{m}_{2}}\mathrm{d}H(u),\]

\noindent
therefore,
\begin{align*}
M_{1}  =&  -\frac{1}{2\pi^{2}}\oint_{\mathcal{C}_{1}}\oint_{\mathcal{C}_{2}}z_{1}^{2}z_{2}^{2}\cdot\frac{1}{1-a}\cdot\frac{\partial^{2}a}{\partial z_{1}\partial z_{2}}\mathrm{d}z_{1}\mathrm{d}z_{2}\\
 = & ~ \frac{2}{2\pi i}\oint_{\mathcal{C}_{2}}z_{2}^{2}\cdot\frac{1}{c}\int_{-1}^{1}\left[\frac{T_{r}(t)T_{s}(t)}{\left(1+T_{s}\left(t\right)\underline{m}_{2}\right)^{2}}\left(-2T_{r}(t)+\frac{1}{c}\int_{-1}^{1}\frac{T_{s}(u)T_{r}(u)\underline{m}_{2}}{1+T_{s}\left(u\right)\underline{m}_{2}}\mathrm{d}H(u)\right)\right]\mathrm{d}H(t)\mathrm{d}\underline{m}_{2}\\
  = & ~\frac{2}{2\pi i}\oint_{\mathcal{C}_{2}}\left[\frac{\left(-1+\frac{1}{c}-\frac{1}{c\sqrt{1-\underline{m}_{2}^{2}}}\right)^{2}}{\underline{m}_{2}^{2}}\cdot\frac{1}{c}\int_{-1}^{1}\frac{-2T_{r}^{2}(t)T_{s}(t)}{\left(1+T_{s}\left(t\right)\underline{m}_{2}\right)^{2}}\mathrm{d}H(t)\right]\mathrm{d}\underline{m}_{2}\\
 + &~ \frac{1}{c}\int_{-1}^{1}\left[\frac{2}{2\pi i}\oint_{\mathcal{C}_{2}}\frac{\left(-1+\frac{1}{c}-\frac{1}{c\sqrt{1-\underline{m}_{2}^{2}}}\right)^{2}T_{r}(t)T_{s}(t)}{\underline{m}_{2}\left(1+T_{s}\left(t\right)\underline{m}_{2}\right)^{2}}\cdot\left(\frac{1}{c}\int_{-1}^{1}\frac{T_{s}(u)T_{r}(u)}{1+T_{s}\left(u\right)\underline{m}_{2}}\mathrm{d}H(u)\right)\mathrm{d}\underline{m}_{2}\right]\mathrm{d}H(t)\\
  = & ~ \frac{1}{c}\int_{-1}^{1}\left(-2T_{r}^{2}(t)T_{s}(t)\right)\cdot\left.\left[\frac{2\left(-1+\frac{1}{c}-\frac{1}{c\sqrt{1-\underline{m}_{2}^{2}}}\right)^{2}}{\left(1+T_{s}\left(t\right)\underline{m}_{2}\right)^{2}}\right]^{(1)}\right\rvert_{\um_2=0}\mathrm{d}H(t)\\
 ~~~&+\left(\frac{1}{c}\int_{-1}^{1}T_{s}(u)T_{r}(u)\mathrm{d}H(u)\right)\cdot\frac{2}{c}\int_{-1}^{1}T_{r}(t)T_{s}(t)\mathrm{d}H(t)\\
  = & ~ \frac{8}{c}\int_{-1}^{1}T_{r}^{2}(t)T_{s}^{2}(t)\mathrm{d}H(t)+\frac{2}{c^2}\left(\int_{-1}^{1}T_{s}(u)T_{r}(u)\mathrm{d}H(u)\right)^{2}.
\end{align*}

\noindent
Similarly, for $M_{2}$, considering
\begin{align*}
~ &  ~ \frac{1}{2\pi i}\oint_{\mathcal{C}_{1}}z_{1}^{2}\cdot\cfrac{1}{\left(1-a\right)^{2}}\cdot\frac{\partial a}{\partial z_{1}}\cdot\frac{\partial a}{\partial z_{2}}\mathrm{d}z_{1}\\
 & =  \frac{1}{2\pi i}\oint_{\mathcal{C}_{1}}\frac{\left(-1+\frac{1}{c}-\frac{1}{c\sqrt{1-\underline{m}_{1}^{2}}}\right)^{2}}{\underline{m}_{1}^{2}\left(1-a\right)^{2}}\cdot\left[\frac{1}{c}\int_{-1}^{1}\frac{T_{s}\left(t\right)\underline{m}_{2}}{1+T_{s}\left(t\right)\underline{m}_{2}}\cdot\frac{T_{r}(t)}{\left(1+T_{r}\left(t\right)\underline{m}_{1}\right)^{2}}\cdot\frac{\partial\underline{m}_{1}}{\partial z_{1}}{\rm \mathrm{d}}H(t)\right]\\
 & ~\cdot  \left[\frac{1}{c}\int_{-1}^{1}\frac{T_{r}\left(t\right)\underline{m}_{1}}{1+T_{r}\left(t\right)\underline{m}_{1}}\cdot\frac{T_{s}(t)}{\left(1+T_{s}\left(t\right)\underline{m}_{2}\right)^{2}}\cdot\frac{\partial\underline{m}_{2}}{\partial z_{2}}{\rm \mathrm{d}}H(t)\right]\mathrm{d}z_{1}\\
 & =  -\frac{\partial\underline{m}_{2}}{\partial z_{2}}\cdot\frac{1}{c}\int_{-1}^{1}\left[\frac{1}{2\pi i}\oint_{\mathcal{C}_{1}}\frac{\left(-1+\frac{1}{c}-\frac{1}{c\sqrt{1-\underline{m}_{1}^{2}}}\right)^{2}}{\underline{m}_{1}\left(1-a\right)^{2}}\cdot\frac{T_{r}(t)T_{s}\left(t\right)\underline{m}_{2}}{\left(1+T_{s}\left(t\right)\underline{m}_{2}\right)\left(1+T_{r}\left(t\right)\underline{m}_{1}\right)^{2}}\right.\\
&~\left. \cdot\left(\frac{1}{c}\int_{-1}^{1}\frac{T_{r}\left(u\right)T_{s}(u)}{\left(1+T_{r}\left(u\right)\underline{m}_{1}\right)\left(1+T_{s}\left(u\right)\underline{m}_{2}\right)^{2}}\mathrm{d}H(u)\right)\mathrm{d}\underline{m}_{1}\right]{\rm \mathrm{d}}H(t)\\
 & =  -\frac{\partial\underline{m}_{2}}{\partial z_{2}}\cdot\frac{1}{c}\int_{-1}^{1}\frac{T_{r}(t)T_{s}\left(t\right)\underline{m}_{2}}{1+T_{s}\left(t\right)\underline{m}_{2}}\mathrm{d}H(t)\cdot\frac{1}{c}\int_{-1}^{1}\frac{T_{r}\left(u\right)T_{s}(u)}{\left(1+T_{s}\left(u\right)\underline{m}_{2}\right)^{2}}\mathrm{d}H(u).
\end{align*}

\noindent
Then,
\begin{align*}
&~M_{2} =  -\frac{1}{2\pi^{2}}\oint_{\mathcal{C}_{1}}\oint_{\mathcal{C}_{2}}z_{1}^{2}z_{2}^{2}\cdot\cfrac{1}{\left(1-a\right)^{2}}\cdot\frac{\partial a}{\partial z_{1}}\cdot\frac{\partial a}{\partial z_{2}}\mathrm{d}z_{1}\mathrm{d}z_{2}\\
 & =  \frac{2}{2\pi i}\oint_{\mathcal{C}_{2}}z_{2}^{2}\left[-\frac{\partial\underline{m}_{2}}{\partial z_{2}}\cdot\frac{1}{c}\int_{-1}^{1}\frac{T_{r}(t)T_{s}\left(t\right)\underline{m}_{2}}{1+T_{s}\left(t\right)\underline{m}_{2}}\mathrm{d}H(t)\cdot\frac{1}{c}\int_{-1}^{1}\frac{T_{r}\left(u\right)T_{s}(u)}{\left(1+T_{s}\left(u\right)\underline{m}_{2}\right)^{2}}\mathrm{d}H(u)\right]\mathrm{d}z_{2}\\
 & =  \frac{1}{c}\int_{-1}^{1}\left[\frac{2}{2\pi i}\oint_{\mathcal{C}_{2}}\frac{T_{r}(t)T_{s}\left(t\right)\left(-1+\frac{1}{c}-\frac{1}{c\sqrt{1-\underline{m}_{2}^{2}}}\right)^{2}}{\underline{m}_{2}\lb 1+T_{s}\left(t\right)\underline{m}_{2}\rb}\left(\frac{1}{c}\int_{-1}^{1}\frac{T_{r}\left(u\right)T_{s}(u)}{\left(1+T_{s}\left(u\right)\underline{m}_{2}\right)^{2}}\mathrm{d}H(u)\right)\mathrm{d}\underline{m}_{2}\right]\mathrm{d}H(t)\\
 & =  \frac{2}{c^{2}}\left(\int_{-1}^{1}T_{r}(t)T_{s}\left(t\right)\mathrm{d}H(t)\right)^{2}.
 \end{align*}

\noindent
As for $M_{3}$,
\begin{align*}
~&~~\frac{1}{2\pi i}\oint_{\mathcal{C}_{1}}z_{1}^{2}\frac{\partial^{2}a}{\partial z_{1}\partial z_{2}}\mathrm{d}z_{1} \\
 &=~~ \frac{1}{2\pi i}\oint_{\mathcal{C}_{1}}\frac{\left(-1+\frac{1}{c}-\frac{1}{c\sqrt{1-\underline{m}_{1}^{2}}}\right)^{2}}{\underline{m}_{1}^{2}}\cdot\left(\frac{1}{c}\int_{-1}^{1}\frac{T_{r}(t)T_{s}(t)}{\left(1+T_{r}\left(t\right)\underline{m}_{1}\right)^{2}\left(1+T_{s}\left(t\right)\underline{m}_{2}\right)^{2}}{\rm \mathrm{d}}H(t)\right)\frac{\partial\underline{m}_{1}}{\partial z_{1}}\cdot\frac{\partial\underline{m}_{2}}{\partial z_{2}}\mathrm{d}z_{1}\\
 & = -\frac{\partial\underline{m}_{2}}{\partial z_{2}}\cdot\frac{1}{c}\int_{-1}^{1}\frac{T_{r}(t)T_{s}(t)}{\left(1+T_{s}\left(t\right)\underline{m}_{2}\right)^{2}}\cdot\left[\frac{1}{2\pi i}\oint_{\mathcal{C}_{1}}\frac{\left(-1+\frac{1}{c}-\frac{1}{c\sqrt{1-\underline{m}_{1}^{2}}}\right)^{2}}{\underline{m}_{1}^{2}\left(1+T_{r}\left(t\right)\underline{m}_{1}\right)^{2}}\mathrm{d}\underline{m}_{1}\right]{\rm \mathrm{d}}H(t)\\
 & =~ -\frac{\partial\underline{m}_{2}}{\partial z_{2}}\cdot\frac{1}{c}\int_{-1}^{1}\frac{T_{r}(t)T_{s}(t)}{\left(1+T_{s}\left(t\right)\underline{m}_{2}\right)^{2}}\cdot\left.\left[\frac{\left(-1+\frac{1}{c}-\frac{1}{c\sqrt{1-\underline{m}_{1}^{2}}}\right)^{2}}{\left(1+T_{r}\left(t\right)\underline{m}_{1}\right)^{2}}\right]^{(1)}\right\rvert_{\underline{m}_{1}=0}{\rm \mathrm{d}}H(t)\\
&  = ~\frac{\partial\underline{m}_{2}}{\partial z_{2}}\cdot\frac{1}{c}\int_{-1}^{1}\frac{2T_{r}^{2}(t)T_{s}(t)}{\left(1+T_{s}\left(t\right)\underline{m}_{2}\right)^{2}}{\rm \mathrm{d}}H(t),
\end{align*}

\noindent
thus
\begin{align*}
M_{3}  = & ~\frac{\beta_{x}}{2\pi i}\oint_{\mathcal{C}_{2}}z_{2}^{2}\left[\frac{\partial\underline{m}_{2}}{\partial z_{2}}\cdot\frac{1}{c}\int_{-1}^{1}\frac{2T_{r}^{2}(t)T_{s}(t)}{\left(1+T_{s}\left(t\right)\underline{m}_{2}\right)^{2}}{\rm \mathrm{d}}H(t)\right]\mathrm{d}z_{2}\\
  = &~ -\frac{\beta_{x}}{c}\int_{-1}^{1}\left[\frac{1}{2\pi i}\oint_{\mathcal{C}_{2}}\frac{\left(-1+\frac{1}{c}-\frac{1}{c\sqrt{1-\underline{m}_{2}^{2}}}\right)^{2}}{\underline{m}_{2}^{2}}\cdot\frac{2T_{r}^{2}(t)T_{s}(t)}{\left(1+T_{s}\left(t\right)\underline{m}_{2}\right)^{2}}\mathrm{d}\underline{m}_{2}\right]{\rm \mathrm{d}}H(t)\\
  = & ~ -\frac{\beta_{x}}{c}\int_{-1}^{1}2T_{r}^{2}(t)T_{s}(t)\left.\left[\frac{\left(-1+\frac{1}{c}-\frac{1}{c\sqrt{1-\underline{m}_{2}^{2}}}\right)^{2}}{\left(1+T_{s}\left(t\right)\underline{m}_{2}\right)^{2}}\right]^{(1)}\right\rvert_{\underline{m}_{2}=0}\mathrm{d}H(t)\\
  = & ~\frac{4\beta_{x}}{c}\int_{-1}^{1}T_{r}^{2}(t)T_{s}^{2}(t)\mathrm{d}H(t).
\end{align*}
Combining the three terms gives 
\[
\Cov\left(X_{f_{r}},X_{f_{s}}\right) 
  = \frac{4(\beta_{x}+2)}{c}\int_{-1}^{1}T_{r}^{2}(t)T_{s}^{2}(t)\mathrm{d}H(t)+\frac{4}{c^{2}}\left(\int_{-1}^{1}T_{s}(u)T_{r}(u)\mathrm{d}H(u)\right)^{2}.
\]
Note that for $\d H(t)=\cfrac{1}{\pi\sqrt{1-t^2}}\d t$, we have 
\begin{align*}
  \int_{-1}^1 T_r(t)T_s(t)\d H(t)
  &=  \frac12 1_{\{r=s\}}, \\
  \int_{-1}^1 T_r^2(t)T_s^2(t)\d H(t)
  & =
    \frac38 1_{\{ r=s\}}  +
    \frac14 1_{\{ r\ne s\}}  .
\end{align*}
Therefore,
\[
\renewcommand{\arraystretch}{1.5}
\Cov\left(X_{f_{r}},X_{f_{s}}\right)=\left\{
  \begin{array}{ll}
    \frac{1+\frac{3}{2}c(\beta_x+2)}{c^2}, &~\mbox{ if } r=s,\\
    \frac{\beta_x+2}{c}, &~\mbox{ if } r\neq s.
  \end{array}
  \right.
\]
The required formula is established.


\begin{table}[htp]
\centering
\footnotesize
\caption{\footnotesize Test Size under Scenario (I) and (II)}\label{Tab:LowSize}
\begin{tabular}{ccc|cccc|cccc}
\hline
\hline
 &  &  & \multicolumn{2}{c}{$\phi_{q}$(I)} & \multicolumn{2}{c|}{$U_{q}$(I)} & \multicolumn{2}{c}{$\phi_{q}$(II)} & \multicolumn{2}{c}{$U_{q}$(II)}\tabularnewline
\cline{4-11}
$p$ & $n$ & $c_{n}$ & $q=1$ & $q=3$ & $q=1$ & $q=3$ & $q=1$ & $q=3$ & $q=1$ & $q=3$\tabularnewline
\hline
5 & 1000 & 0.005 & 0.081 & 0.078 & 0.065 & 0.049 & 0.081 & 0.074 & 0.090 & 0.081\tabularnewline
10 & 2000 & 0.005 & 0.059 & 0.060 & 0.052 & 0.050 & 0.058 & 0.055 & 0.068 & 0.067\tabularnewline
25 & 5000 & 0.005 & 0.051 & 0.054 & 0.047 & 0.043 & 0.054 & 0.059 & 0.055 & 0.050\tabularnewline
40 & 8000 & 0.005 & 0.050 & 0.049 & 0.054 & 0.036 & 0.062 & 0.055 & 0.057 & 0.051\tabularnewline
\hline
10 & 1000 & 0.01 & 0.072 & 0.067 & 0.057 & 0.048 & 0.063 & 0.060 & 0.068 & 0.052\tabularnewline
20 & 2000 & 0.01 & 0.066 & 0.059 & 0.050 & 0.040 & 0.057 & 0.058 & 0.065 & 0.049\tabularnewline
50 & 5000 & 0.01 & 0.046 & 0.053 & 0.045 & 0.044 & 0.048 & 0.053 & 0.055 & 0.044\tabularnewline
80 & 8000 & 0.01 & 0.056 & 0.049 & 0.051 & 0.046 & 0.046 & 0.046 & 0.052 & 0.050\tabularnewline
\hline
50 & 1000 & 0.05 & 0.063 & 0.056 & 0.051 & 0.048 & 0.046 & 0.058 & 0.055 & 0.058\tabularnewline
100 & 2000 & 0.05 & 0.058 & 0.052 & 0.061 & 0.047 & 0.048 & 0.052 & 0.049 & 0.047\tabularnewline
250 & 5000 & 0.05 & 0.056 & 0.055 & 0.050 & 0.047 & 0.046 & 0.048 & 0.047 & 0.037\tabularnewline
400 & 8000 & 0.05 & 0.049 & 0.047 & 0.047 & 0.040 & 0.055 & 0.042 & 0.046 & 0.044\tabularnewline

\hline
10 & 100 & 0.1 & 0.073 & 0.075 & 0.062 & 0.061 & 0.072 & 0.074 & 0.079 & 0.083\tabularnewline
40 & 400 & 0.1 & 0.053 & 0.062 & 0.050 & 0.043 & 0.051 & 0.059 & 0.056 & 0.055\tabularnewline
60 & 600 & 0.1 & 0.049 & 0.043 & 0.055 & 0.047 & 0.052 & 0.053 & 0.053 & 0.049\tabularnewline
100 & 1000 & 0.1 & 0.062 & 0.058 & 0.053 & 0.045 & 0.051 & 0.054 & 0.054 & 0.044\tabularnewline
\hline
50 & 100 & 0.5 & 0.066 & 0.067 & 0.066 & 0.050 & 0.051 & 0.059 & 0.069 & 0.070\tabularnewline
200 & 400 & 0.5 & 0.053 & 0.052 & 0.051 & 0.045 & 0.048 & 0.045 & 0.053 & 0.059\tabularnewline
300 & 600 & 0.5 & 0.053 & 0.054 & 0.044 & 0.046 & 0.040 & 0.048 & 0.043 & 0.038\tabularnewline
500 & 1000 & 0.5 & 0.052 & 0.050 & 0.045 & 0.051 & 0.041 & 0.045 & 0.048 & 0.038\tabularnewline
\hline
90 & 100 & 0.9 & 0.051 & 0.055 & 0.050 & 0.057 & 0.048 & 0.051 & 0.069 & 0.078\tabularnewline
360 & 400 & 0.9 & 0.056 & 0.051 & 0.050 & 0.050 & 0.047 & 0.048 & 0.047 & 0.047\tabularnewline
540 & 600 & 0.9 & 0.058 & 0.050 & 0.061 & 0.049 & 0.037 & 0.046 & 0.054 & 0.057\tabularnewline
900 & 1000 & 0.9 & 0.048 & 0.053 & 0.058 & 0.046 & 0.045 & 0.045 & 0.048 & 0.045\tabularnewline
\hline
150 & 100 & 1.5 & 0.042 & 0.047 & 0.065 & 0.064 & 0.039 & 0.048 & 0.059 & 0.070\tabularnewline
600 & 400 & 1.5 & 0.048 & 0.054 & 0.049 & 0.040 & 0.036 & 0.045 & 0.048 & 0.049\tabularnewline
900 & 600 & 1.5 & 0.056 & 0.055 & 0.040 & 0.046 & 0.039 & 0.043 & 0.048 & 0.046\tabularnewline
1500 & 1000 & 1.5 & 0.051 & 0.052 & 0.041 & 0.043 & 0.042 & 0.049 & 0.049 & 0.047\tabularnewline
\hline
200 & 100 & 2 & 0.051 & 0.051 & 0.057 & 0.049 & 0.045 & 0.045 & 0.067 & 0.066\tabularnewline
800 & 400 & 2 & 0.047 & 0.056 & 0.052 & 0.049 & 0.043 & 0.052 & 0.046 & 0.048\tabularnewline
1200 & 600 & 2 & 0.055 & 0.050 & 0.051 & 0.052 & 0.043 & 0.052 & 0.053 & 0.040\tabularnewline
2000 & 1000 & 2 & 0.054 & 0.053 & 0.050 & 0.051 & 0.048 & 0.048 & 0.049 & 0.040\tabularnewline
\hline
500 & 100 & 5 & 0.063 & 0.053 & 0.049 & 0.044 & 0.051 & 0.052 & 0.051 & 0.068\tabularnewline
2000 & 400 & 5 & 0.049 & 0.054 & 0.042 & 0.044 & 0.040 & 0.048 & 0.055 & 0.048\tabularnewline
3000 & 600 & 5 & 0.052 & 0.056 & 0.052 & 0.042 & 0.034 & 0.044 & 0.049 & 0.051\tabularnewline
5000 & 1000 & 5 & 0.052 & 0.052 & 0.047 & 0.049 & 0.043 & 0.045 & 0.049 & 0.057\tabularnewline
\hline
\hline
\end{tabular}
\end{table}

\begin{table}[h]
\centering
\footnotesize
\caption{\footnotesize Test Power under Scenario (III) and (IV)}\label{Tab:LowPower}
\begin{tabular}{ccc|cccc|cccc}
\hline
\hline
 &  &  & \multicolumn{2}{c}{$\phi_{q}$(III)} & \multicolumn{2}{c|}{$U_{q}$(III)} & \multicolumn{2}{c}{$\phi_{q}$(IV)} & \multicolumn{2}{c}{$U_{q}$(IV)}\tabularnewline
\cline{4-11}
$p$ & $n$ & $c_{n}$ & $q=1$ & $q=3$ & $q=1$ & $q=3$ & $q=1$ & $q=3$ & $q=1$ & $q=3$\tabularnewline
\hline
5 & 1000 & 0.005 & 0.999 & 0.990 & 0.807 & 0.635 & 0.999 & 0.986 & 0.761 & 0.633\tabularnewline
10 & 2000 & 0.005 & 1 & 1 & 0.9995 & 0.986 & 1 & 1 & 0.999 & 0.984\tabularnewline
25 & 5000 & 0.005 & 1 & 1 & 1 & 1 & 1 & 1 & 1 & 1\tabularnewline
40 & 8000 & 0.005 & 1 & 1 & 1 & 1 & 1 & 1 & 1 & 1\tabularnewline
\hline
10 & 1000 & 0.01 & 1 & 0.998 & 0.824 & 0.647 & 1 & 0.997 & 0.79 & 0.636\tabularnewline
20 & 2000 & 0.01 & 1 & 1 & 1 & 0.992 & 1 & 1 & 1 & 0.989\tabularnewline
50 & 5000 & 0.01 & 1 & 1 & 1 & 1 & 1 & 1 & 1 & 1\tabularnewline
80 & 8000 & 0.01 & 1 & 1 & 1 & 1 & 1 & 1 & 1 & 1\tabularnewline
\hline
50 & 1000 & 0.05 & 1 & 0.9995 & 0.859 & 0.681 & 1 & 1 & 0.835 & 0.660\tabularnewline
100 & 2000 & 0.05 & 1 & 1 & 1 & 0.993 & 1 & 1 & 1 & 0.996\tabularnewline
250 & 5000 & 0.05 & 1 & 1 & 1 & 1 & 1 & 1 & 1 & 1\tabularnewline
400 & 8000 & 0.05 & 1 & 1 & 1 & 1 & 1 & 1 & 1 & 1\tabularnewline
\hline

10 & 100 & 0.1 & 0.335 & 0.213 & 0.085 & 0.077 & 0.3055 & 0.1955 & 0.129 & 0.128\tabularnewline
40 & 400 & 0.1 & 0.983 & 0.764 & 0.294 & 0.207 & 0.968 & 0.6935 & 0.283 & 0.225\tabularnewline
60 & 600 & 0.1 & 1 & 0.973 & 0.482 & 0.338 & 1 & 0.94 & 0.483 & 0.331\tabularnewline
100 & 1000 & 0.1 & 1 & 1 & 0.851 & 0.651 & 1 & 1 & 0.852 & 0.652\tabularnewline
\hline
50 & 100 & 0.5 & 0.416 & 0.245 & 0.093 & 0.093 & 0.3075 & 0.1805 & 0.094 & 0.091\tabularnewline
200 & 400 & 0.5 & 1 & 0.957 & 0.287 & 0.190 & 0.9995 & 0.8585 & 0.310 & 0.223\tabularnewline
300 & 600 & 0.5 & 1 & 1 & 0.497 & 0.327 & 1 & 0.995 & 0.491 & 0.322\tabularnewline
500 & 1000 & 0.5 & 1 & 1 & 0.878 & 0.669 & 1 & 1 & 0.886 & 0.646\tabularnewline
\hline
90 & 100 & 0.9 & 0.506 & 0.311 & 0.090 & 0.079 & 0.3785 & 0.217 & 0.096 & 0.097\tabularnewline
360 & 400 & 0.9 & 1 & 0.991 & 0.323 & 0.208 & 1 & 0.952 & 0.315 & 0.214\tabularnewline
540 & 600 & 0.9 & 1 & 1 & 0.527 & 0.338 & 1 & 0.9995 & 0.504 & 0.331\tabularnewline
900 & 1000 & 0.9 & 1 & 1 & 0.897 & 0.655 & 1 & 1 & 0.894 & 0.671\tabularnewline
\hline
150 & 100 & 1.5 & 0.607 & 0.365 & 0.089 & 0.078 & 0.4545 & 0.2745 & 0.102 & 0.098\tabularnewline
600 & 400 & 1.5 & 1 & 1 & 0.337 & 0.203 & 1 & 0.992 & 0.319 & 0.212\tabularnewline
900 & 600 & 1.5 & 1 & 1 & 0.573 & 0.346 & 1 & 1 & 0.576 & 0.326\tabularnewline
1500 & 1000 & 1.5 & 1 & 1 & 0.922 & 0.674 & 1 & 1 & 0.906 & 0.656\tabularnewline
\hline
200 & 100 & 2 & 0.694 & 0.428 & 0.092 & 0.082 & 0.5165 & 0.2985 & 0.108 & 0.098\tabularnewline
800 & 400 & 2 & 1 & 1 & 0.351 & 0.207 & 1 & 0.9985 & 0.352 & 0.204\tabularnewline
1200 & 600 & 2 & 1 & 1 & 0.609 & 0.351 & 1 & 1 & 0.592 & 0.342\tabularnewline
2000 & 1000 & 2 & 1 & 1 & 0.923 & 0.657 & 1 & 1 & 0.933 & 0.656\tabularnewline
\hline
500 & 100 & 5 & 0.9375 & 0.704 & 0.102 & 0.082 & 0.809 & 0.519 & 0.108 & 0.097\tabularnewline
2000 & 400 & 5 & 1 & 1 & 0.452 & 0.190 & 1 & 1 & 0.430 & 0.202\tabularnewline
3000 & 600 & 5 & 1 & 1 & 0.734 & 0.353 & 1 & 1 & 0.712 & 0.329\tabularnewline
5000 & 1000 & 5 & 1 & 1 & 0.986 & 0.659 & 1 & 1 & 0.984 & 0.679\tabularnewline
\hline
\hline
\end{tabular}
\end{table}

\renewcommand{\arraystretch}{1.4}
\begin{table}[!h]
{
\caption{Size and power comparison of permutation test and $\phi_q$ for Gaussian distribution under Scenario (III), nominal level $\alpha=5\%$, testing size is shown in the first row of each $(p,n)$ configuration block. }\label{tab:compareGauss}
\begin{tabular}{cccc|cc|cc}
	\hline 
	&  &  &  & \multicolumn{2}{c|}{Permutation} & \multicolumn{2}{c}{$\phi_{q}$}\tabularnewline
	\cline{5-8} 
	$p$ & $n$ & $c_{n}$ & $a$ & $q=1$ & $q=3$ & $q=1$ & $q=3$\tabularnewline
	\hline 
	150 & 300 & 0.5 & {\bf 0} & {\bf 0.056} & {\bf 0.062} & {\bf 0.052} & {\bf 0.054}\tabularnewline
	150 & 300 & 0.5 & 0.05 & 0.360 & 0.256 & 0.254 & 0.146\tabularnewline
	150 & 300 & 0.5 & 0.09 & 0.992 & 0.948 & 0.934 & 0.694\tabularnewline
	150 & 300 & 0.5 & 0.1 & 1 & 0.988 & 0.976 & 0.820\tabularnewline
	\hline 
	270 & 300 & 0.9 & {\bf 0} & {\bf 0.068} & {\bf 0.074} & {\bf 0.060} & {\bf 0.050}\tabularnewline
	270 & 300 & 0.9 & 0.05 & 0.510 & 0.408 & 0.376 & 0.216\tabularnewline
	270 & 300 & 0.9 & 0.09 & 1 & 1 & 0.990 & 0.820\tabularnewline
	270 & 300 & 0.9 & 0.1 & 1 & 1 & 0.998 & 0.914\tabularnewline
	\hline 
	600 & 300 & 2 & {\bf 0} & {\bf 0.056} & {\bf 0.066} & {\bf 0.062} & {\bf 0.062}\tabularnewline
	600 & 300 & 2 & 0.05 & 0.808 & 0.708 & 0.536 & 0.294\tabularnewline
	600 & 300 & 2 & 0.09 & 1 & 1 & 1 & 0.974\tabularnewline
	600 & 300 & 2 & 0.1 & 1 & 1 & 1 & 0.994\tabularnewline
	\hline 
\end{tabular}
}
\end{table}

\renewcommand{\arraystretch}{1.4}
\begin{table}[!h]
{
\caption{Size and power comparison of permutation test and $\phi_q$ for Gamma distribution under Scenario (IV), nominal level $\alpha=5\%$, testing size is shown in the first row of each $(p,n)$ configuration block. }\label{tab:compareNonGauss}
\begin{tabular}{cccc|cc|cc}
	\hline 
	&  &  &  & \multicolumn{2}{c|}{Permutation} & \multicolumn{2}{c}{$\phi_{q}$}\tabularnewline
	\cline{5-8} 
	$p$ & $n$ & $c_{n}$ & $a$ & $q=1$ & $q=3$ & $q=1$ & $q=3$\tabularnewline
	\hline 
	150 & 300 & 0.5 & {\bf 0} & {\bf 0.052} & {\bf 0.044} & {\bf 0.050} & {\bf 0.058}\tabularnewline
	150 & 300 & 0.5 & 0.05 & 0.344 & 0.276 & 0.228 & 0.142\tabularnewline
	150 & 300 & 0.5 & 0.09 & 0.996 & 0.962 & 0.892 & 0.534\tabularnewline
	150 & 300 & 0.5 & 0.1 & 1 & 0.988 & 0.956 & 0.692\tabularnewline
	\hline 
	270 & 300 & 0.9 & {\bf 0} & {\bf 0.058} & {\bf 0.048} & {\bf 0.046} & {\bf 0.052}\tabularnewline
	270 & 300 & 0.9 & 0.05 & 0.470 & 0.372 & 0.244 & 0.144\tabularnewline
	270 & 300 & 0.9 & 0.09 & 1 & 1 & 0.960 & 0.662\tabularnewline
	270 & 300 & 0.9 & 0.1 & 1 & 1 & 0.988 & 0.802\tabularnewline
	\hline 
	600 & 300 & 2 & {\bf 0} & {\bf 0.058} & {\bf 0.054} & {\bf 0.042} & {\bf 0.056}\tabularnewline
	600 & 300 & 2 & 0.05 & 0.766 & 0.692 & 0.366 & 0.202\tabularnewline
	600 & 300 & 2 & 0.09 & 1 & 1 & 0.998 & 0.884\tabularnewline
	600 & 300 & 2 & 0.1 & 1 & 1 & 1 & 0.944\tabularnewline
	\hline 
\end{tabular}
}
\end{table}


\begin{thebibliography}{99}
\bibitem[{Anderson(2003)}]{Anderson2003} Anderson, T. W.  (2003).  {\em An Introduction to Multivariate
    Statistical Analysis}.  3rd Edition.  John Wiley \& Sons.




\bibitem[{Bai(1999)}]{bai1999}
  {Bai, Z. D.}, (1999) Methodologies in spectral analysis of large dimensional random matrices, a review. \textit{Statistica Sinica}, {\bf 9}, 611-677.



\bibitem[{Bai et al.(2009)}]{bjyz09}
  {Bai, Z. D.}, {Jiang, D. D.}, {Yao, J. F.} and {Zheng, S. R.} (2009).
  Corrections to LRT on large dimensional covariance matrix by RMT. \textit{Annals of Statistics}, {\bf 37}, 3822-3840.

\bibitem[{Bai et al.(2013)}]{bjyz13}
  {Bai, Z. D.}, {Jiang, D. D.}, Yao, J. F., and Zheng, S. R. (2013).
  Testing linear hypotheses in high-dimensional regressions. \textit{Statistics}, {\bf 47}, 1207-1223.

\bibitem[{Bai and Silverstein(1998)}]{BS98}{Bai, Z. D.} and {Silverstein, J. W.} (1998). No eigenvalues outside the support of the limiting
  spectral distribution of large-dimensional sample covariance matrices. \textit{Ann. Probab.}, {\bf 26}, 316-345.

\bibitem[{Bai and Silverstein(2004)}]{BS04}{Bai, Z. D.} and {Silverstein, J. W.} (2004). CLT for linear spectral
  statistics of large dimensional sample covariance matrices. \textit{Ann. Probab.}, {\bf 32}, 553-605.

\bibitem[{Bai and Silverstein(2010)}]{BS10}
  {Bai, Z. D.} and {Silverstein, J. W.} (2010). {\sl Spectral Analysis
    of Large Dimensional Random Matrices}. {Science Press: Beijing}.

  \bibitem[Bai and Wang(2015)]{Bai15}
{Bai, Z. D. and Wang, C.} (2015). A note on the limiting spectral distribution of a symmetrized auto-cross covariance matrix.  {\em Statistics \& Probability Letters}, {\bf 96}, 333-340.


\bibitem[{Bai and Zhou(2008)}]{BZ2008} {Bai, Z. D.} and {Zhou, W.} (2008). Large sample covariance matrices without independence structures in columns. \textit{Statistica Sinica}, {\bf 28}, 425-442.


\bibitem[Chen and Pan(2015)]{CP15}
Chen, B. B. and Pan, G. M. (2015)  CLT for linear spectral statistics of normalized sample covariance matrices with the dimension much larger than the sample size.
\textit{Bernoulli}, {\bf  21}, 1089-1133.


\bibitem[{Billingsley(1995)}]{Billingsley1995} {Billingsley, P.} (1995). \textit{Probability and Measure}. John Wiley $\&$ Sons: New York

\bibitem[{Burkholder(1973)}]{Burkholder1973} {Burkholder, D. L.} (1973). Distribution function inequalities for martingales. \textit{ Ann. Probab.}, {\bf 1}, 19-42.

\bibitem[Gray(2006)]{Gray06}
{Gray, R. M. }(2006). Toeplitz and Circulant Matrices: a Review. {\em Now Publishers Inc.}


\bibitem[Grenander and Szeg\"{o}(1958)]{Grenander58}
{Grenander, U. and Szeg\"{o}, G.}(1958). Toeplitz Forms and Their Applications. In: California Monographs in Mathematical Sciences. {\em University of California Press}, Berkeley.

\bibitem[Johnstone(2007)]{John07}
  {Johnstone, I. M.} (2007). High dimensional statistical inference and random matrices.
  \textit{Int. Cong. Mathematicians}, \textbf{Vol. I}, 307-333. Z$\ddot{u}$rich, Switzerland: European Mathematical
  Society.

\bibitem[{Li et al.(2016)}]{Li16}
{Li, Z., Yao, J., Lam, C., and Yao, Q. }(2016). On testing a high-dimensional white noise. \textit{Manuscript}.

\bibitem[{Mar\v{c}enko and Pastur(1967)}]{MP}
  {Mar\v{c}enko, V. A.} and {Pastur, L. A.} (1967).
  \newblock Distribution of eigenvalues for some sets of random matrices.
  \newblock {\em Math. USSR-Sb}, {\bf 1}, 457-483.


\bibitem[Pan(2012)]{Pan12}
  {Pan, G. M.} (2012).
  Comparison between two types of large sample covariance matrices. \textit{Annales de l'Institut Henri Poincare-Probabiliteset Statistiques}, {\bf 50}, 655-677.
  
  \bibitem[Paul and Aue(2014)]{PA14}
  Paul, D. and Aue, A. (2014). Random matrix theory in statistics: A review, \textit{Journal of Statistical Planning and Inference}, {\bf 150}, 1-29.

\bibitem[{Silverstein(1995)}]{Silv95} {Silverstein,  J. W.} (1995).  Strong convergence of the empirical
  distribution of eigenvalues of large dimensional random matrices. \textit{J. Multivariate Anal.}, {\bf 5},  331-339.

\bibitem[{Silverstein and Bai(1995)}]{SB95}{Silverstein, J. W.} and {Bai, Z. D.} (1995). On the empirical distribution of eigenvalues
  of a class of large dimensional random matrices. \textit{J. Multivariate Anal.}, {\bf 54}, 175-192.

\bibitem[{Silverstein and Choi(1995)}]{SC1995} {Silverstein, J. W.} and {Choi, S. I.} (1995).
  Analysis of the limiting spectral distribution of large dimensional random matrices.
  \textit{J. Multivariate Anal.}, {\bf 54}, 295-309.


\bibitem[Simes(1986)]{Simes86}
\rdd{Simes, R. J. }(1986). An improved Bonferroni procedure for multiple tests of significance. {\em Biometrika}, {\bf73}, 751-754.

\bibitem[{Srivastava(2005)}]{Sriv05}
  \mgai{
  {Srivastava, M. S.} (2005).
  \newblock Some tests concerning the covariance matrix in high dimensional data.
  \newblock \textit{J. Japan Statist. Soc.}, {\bf 35}, 251--272.}

\bibitem[Yao et al.(2015)]{Yao15}
Yao, J. F., Bai, Z. D., and Zheng, S. R. (2015). {\sl Large Sample Covariance Matrices and High-Dimensional Data Analysis}, Cambridge University Press.

\bibitem[{Zheng et al.(2015)}]{ZBY2015}
  {Zheng, S. R.}, {Bai, Z. D.}, and {Yao, J. F.} (2015). Substitution principle for CLT of linear spectral statistics
  of high-dimensional sample covariance matrices with applications to hypothesis testing. \textit{Ann.  Statist.}, {\bf 43}, 546-591.


\bibitem[{Zheng et al.(2017a)}]{ZBY-cltF}
  {Zheng, S. R.}, {Bai, Z. D.}, and {Yao, J. F.} (2017a). CLT for large dimensional general Fisher matrices and its applications
  in high-dimensional data analysis. \textit{Bernoulli},
  {\bf 23}, 1130-1178.

 \bibitem[Zheng et al.(2017b)]{ZBYZ16}
 {Zheng, S. R.}, {Bai, Z. D.}, Yao, J. F., and Zhu, H. T. (2017b). CLT for linear spectral statistics of large dimensional sample covariance matrices with dependent data. \textit{Manuscript}.



\end{thebibliography}
\end{document}